\documentclass[12pt]{article}
\usepackage{latexsym}
\usepackage{amsfonts}
\parskip=\medskipamount                 
\thicklines                     
\newcommand{\bimn}[7]{\bibitem{#1}#2,
{\em #3},
{ #4}\hspace{0.25em}{\bf
#5}\hspace{0.25em}(#6)\hspace{0.25em}{#7}.}

\def\inbar{\vrule height1.5ex width.4pt depth0pt}
\def\IC{\relax\,\hbox{$\inbar\kern-.3em{\rm C}$}}
\def\IN{\relax{\rm I\kern-.18em N}}
\def\IQ{\relax\,\hbox{$\inbar\kern-.3em{\rm Q}$}}
\def\IR{\relax{\rm I\kern-.18em R}}
\def\ZZ{\relax{\sf Z\kern-.4em Z}}
\def\a{\alpha} \def\b{\beta}    
 \def\l{\lambda} 
  \def\s{\sigma}

 \def\cK{{\cal K}} \def\cL{{\cal L}} \def\cM{{\cal M}}

\newtheorem{theorem}{Theorem}[section]
\newtheorem{proposition}[theorem]{Proposition}
\newtheorem{corollary}[theorem]{Corollary}
\newtheorem{conjecture}[theorem]{Conjecture}
\newtheorem{lemma}[theorem]{Lemma}
\newtheorem{definition}[theorem]{Definition}
\newtheorem{remark}[theorem]{Remark}
\newtheorem{mt}{Main Theorem}

\catcode`\@=11

\newif\if@fewtab\@fewtabtrue

\catcode`\@=11

\newif\if@fewtab\@fewtabtrue

{\count255=\time\divide\count255 by 60
\xdef\hourmin{\number\count255}
\multiply\count255 by-60\advance\count255 by\time
\xdef\hourmin{\hourmin:\ifnum\count255<10 0\fi\the\count255}}
\def\ps@draft{\let\@mkboth\@gobbletwo
    \def\@oddhead{}
    \def\@oddfoot
      {\hbox to 7 cm{\footnotesize {\em Draft of \jobname:} \draftdate
       \hfil}\hskip -7cm\hfil\rm\thepage \hfil}
    \def\@evenhead{}\let\@evenfoot\@oddfoot}

\def\ceqno{\global\@fewtabfalse
    \ifcase\@eqcnt \def\@tempa{& & &}\or \def\@tempa{& &}
      \or \def\@tempa{&}
      \or\def\@tempa{}\fi\@tempa
{\rm(\theequation)}}

\def\aeqno#1{\global\@fewtabfalse
    \ifcase\@eqcnt \def\@tempa{& & &}\or \def\@tempa{& &}
      \or \def\@tempa{&}
      \or\def\@tempa{}\fi\@tempa
{\rm(\theequation,#1)}}

\def\label#1{\ifnum\draftcontrol=1
 \global\def\draftnote{$\scriptstyle #1$}\fi
 \@bsphack\if@filesw {\let\thepage\relax
   \def\protect{\noexpand\noexpand\noexpand}%
\xdef\@gtempa{\write\@auxout{\string
      \newlabel{#1}{{\@currentlabel}{\thepage}}}}}\@gtempa
   \if@nobreak \ifvmode\nobreak\fi\fi\fi
  \@esphack}

\def\alabel#1#2{\label{#1}\global\@fewtabfalse
    \ifcase\@eqcnt \def\@tempa{& & &}\or \def\@tempa{& &}
      \or \def\@tempa{&}
      \or\def\@tempa{}\fi\@tempa
{\hbox to 3cm{\phantom{\rm(\theequation,#2)}
\draftnote \hfil}\hskip -3cm {\rm(\theequation,#2)}}}

\def\clabel#1{\label{#1}\global\@fewtabfalse
    \ifcase\@eqcnt \def\@tempa{& & &}\or \def\@tempa{& &}
      \or \def\@tempa{&}
      \or\def\@tempa{}\fi\@tempa
{\hbox to 3cm{\phantom{\rm(\theequation)}
\draftnote \hfil}\hskip -3cm{\rm(\theequation)}}}

\def\eqnarray{\def\draftnote{{}}\global\@fewtabtrue
\stepcounter{equation}\let\@currentlabel=\theequation
\global\@eqnswtrue
\global\@eqcnt\z@\tabskip\@centering\let\\=\@eqncr
$$\halign to \displaywidth\bgroup\@eqnsel\hskip\@centering\@eqcnt\z@
  $\displaystyle\tabskip\z@{##}$&\global\@eqcnt\@ne
  \hskip 1\arraycolsep \hfil$\displaystyle{##}$\hfil
  &\global\@eqcnt\tw@ \hskip 1\arraycolsep
$\displaystyle\tabskip\z@{##}$
\hfil  \tabskip\@centering&\global\@eqcnt\thr@@\llap{##}\tabskip\z@
\cr}

\def\endeqnarray{\@@eqncr\egroup
      \global\advance\c@equation\m@ne$$\global\@ignoretrue}

\def\@eqnnum{\hbox to 3cm{\phantom{\rm(\theequation)} \draftnote
                         \hfil}\hskip -3cm {\rm(\theequation)}}

\def\@@eqncr{\let\@tempa\relax
    \ifcase\@eqcnt \def\@tempa{& & &}\or \def\@tempa{& &}
      \or \def\@tempa{&}
      \or\def\@tempa{}
\fi\@tempa
\if@eqnsw
\if@fewtab\@eqnnum\fi
\stepcounter{equation}\fi\global
\@eqnswtrue\global\@eqcnt\z@\global\@fewtabtrue\cr}

\def\draftcite#1{\ifnum\draftcontrol=1#1\else{}\fi}

\def\@lbibitem[#1]#2{\item{}\hskip -3cm \hbox to 2cm
{\hfil$\scriptstyle\draftcite{#2}$}\hskip
1cm[\@biblabel{#1}]\if@filesw
     {\def\protect##1{\string ##1\space}\immediate
      \write\@auxout{\string\bibcite{#2}{#1}}}\fi\ignorespaces}

\def\@bibitem#1{\item\hskip -3cm \hbox to 2cm
{\hfil $\scriptstyle\draftcite{#1}$}\hskip 1cm
\if@filesw \immediate\write\@auxout
       {\string\bibcite{#1}{\the\value{\@listctr}}}\fi\ignorespaces}

\def\nsection#1{\section{#1}\setcounter{equation}{0}}

\def\draftdate{\number\month/\number\day/\number\year\ \ \ \hourmin }

\global\def\draftcontrol{0}
\catcode`\@=12

\def\theequation{{\thesection.\arabic{equation}}}

\def\qq{\begin{eqnarray}}
\def\qqq{\end{eqnarray}}
\def\rx#1{~(\ref{#1})}
\def\ex#1{eq.\hspace*{-3pt}\rx{#1}}
\def\eex#1{eqs.\hspace*{-3pt}\rx{#1}}
\def\cx#1{~\cite{#1}}
\def\rw#1{~\ref{#1}}

\newlength{\shiftwidth}
\def\shift#1{&&\hbox to \shiftwidth{\hfill $\displaystyle#1$}}
\newlength{\sshiftwidth}
\def\sshift#1{\lefteqn{\hbox to
\sshiftwidth{\hfill$\displaystyle#1$}}}

\def\ie{{\it i.e.\ }}

\def\cf{{\it cf.\ }}
\def\rhs{{\it r.h.s.\ }}
\def\lhs{{\it l.h.s.\ }}

\def\p{^{\prime}}

\def\prosign{\mathop{{\rm sign}}\nolimits}

\def\promod{\mathop{{\rm mod}}\nolimits}
\def\mod#1{\;(\promod #1)}

\def\pr#1#2{ \noindent{\em Proof of #1~\ref{#2}.} }
\def\proof{ \noindent{\em Proof.} }

\def\qed{ \hfill $\Box$ }

\def\lrbc#1{ \left( #1 \right) }
\def\lrbs#1{ \left[ #1 \right] }

\def\u#1{ \underline{#1} }
\def\uu#1{ \underline{\underline{#1}} }

\def\ux{ {\u{x}} }
\def\uy{ {\u{y}} }

\def\ual{ {\u{\a}} }

\def\um{ {\u{m}} }

\def\ZZ{ \mathbb{Z} }
\def\IQ{ \mathbb{Q} }
\def\IC{ \mathbb{C} }

\def\qbezier{\bezier{120}}
\def\DottedCircle{
\bezier{4}(0.966,-0.259)(1.04,0)(0.966,0.259)
\bezier{4}(0.966,0.259)(0.897,0.518)(0.707,0.707)
\bezier{4}(0.707,0.707)(0.518,0.897)(0.259,0.966)
\bezier{4}(0.259,0.966)(0,1.04)(-0.259,0.966)
\bezier{4}(-0.259,0.966)(-0.518,0.897)(-0.707,0.707)
\bezier{4}(-0.707,0.707)(-0.897,0.518)(-0.966,0.259)
\bezier{4}(-0.966,0.259)(-1.04,0)(-0.966,-0.259)
\bezier{4}(-0.966,-0.259)(-0.897,-0.518)(-0.707,-0.707)
\bezier{4}(-0.707,-0.707)(-0.518,-0.897)(-0.259,-0.966)
\bezier{4}(-0.259,-0.966)(0,-1.04)(0.259,-0.966)
\bezier{4}(0.259,-0.966)(0.518,-0.897)(0.707,-0.707)
\bezier{4}(0.707,-0.707)(0.897,-0.518)(0.966,-0.259)
}
\def\Endpoint[#1]{
\ifcase#1
\put(1,0){\circle*{0.15}}
\or\put(0.866,0.5){\circle*{0.15}}
\or\put(0.5,0.866){\circle*{0.15}}
\or\put(0,1){\circle*{0.15}}
\or\put(-0.5,0.866){\circle*{0.15}}
\or\put(-0.866,0.5){\circle*{0.15}}
\or\put(-1,0){\circle*{0.15}}
\or\put(-0.866,-0.5){\circle*{0.15}}
\or\put(-0.5,-0.866){\circle*{0.15}}
\or\put(0,-1){\circle*{0.15}}
\or\put(0.5,-0.866){\circle*{0.15}}
\or\put(0.866,-0.5){\circle*{0.15}}
\fi}
\def\Arc[#1]{
\thicklines			
\ifcase#1
\bezier{25}(0.966,-0.259)(1.04,0)(0.966,0.259)
\or
\bezier{25}(0.966,0.259)(0.897,0.518)(0.707,0.707)
\or
\bezier{25}(0.707,0.707)(0.518,0.897)(0.259,0.966)
\or
\bezier{25}(0.259,0.966)(0,1.04)(-0.259,0.966)
\or
\bezier{25}(-0.259,0.966)(-0.518,0.897)(-0.707,0.707)
\or
\bezier{25}(-0.707,0.707)(-0.897,0.518)(-0.966,0.259)
\or
\bezier{25}(-0.966,0.259)(-1.04,0)(-0.966,-0.259)
\or
\bezier{25}(-0.966,-0.259)(-0.897,-0.518)(-0.707,-0.707)
\or
\bezier{25}(-0.707,-0.707)(-0.518,-0.897)(-0.259,-0.966)
\or
\bezier{25}(-0.259,-0.966)(0,-1.04)(0.259,-0.966)
\or
\bezier{25}(0.259,-0.966)(0.518,-0.897)(0.707,-0.707)
\or
\bezier{25}(0.707,-0.707)(0.897,-0.518)(0.966,-0.259)
\fi}
\def\DottedArc[#1]{
\ifcase#1
\bezier{4}(0.966,-0.259)(1.04,0)(0.966,0.259)
\or
\bezier{4}(0.966,0.259)(0.897,0.518)(0.707,0.707)
\or
\bezier{4}(0.707,0.707)(0.518,0.897)(0.259,0.966)
\or
\bezier{4}(0.259,0.966)(0,1.04)(-0.259,0.966)
\or
\bezier{4}(-0.259,0.966)(-0.518,0.897)(-0.707,0.707)
\or
\bezier{4}(-0.707,0.707)(-0.897,0.518)(-0.966,0.259)
\or
\bezier{4}(-0.966,0.259)(-1.04,0)(-0.966,-0.259)
\or
\bezier{4}(-0.966,-0.259)(-0.897,-0.518)(-0.707,-0.707)
\or
\bezier{4}(-0.707,-0.707)(-0.518,-0.897)(-0.259,-0.966)
\or
\bezier{4}(-0.259,-0.966)(0,-1.04)(0.259,-0.966)
\or
\bezier{4}(0.259,-0.966)(0.518,-0.897)(0.707,-0.707)
\or
\bezier{4}(0.707,-0.707)(0.897,-0.518)(0.966,-0.259)
\fi}
\def\Chord[#1,#2]{
\thinlines
\ifnum#1>#2\Chord[#2,#1]
\else\ifnum#1<#2
\ifcase#1
\ifcase#2
\or\qbezier(1,0)(0.516,0.138)(0.866,0.5)
\or\qbezier(1,0)(0.45,0.26)(0.5,0.866)
\or\qbezier(1,0)(0.327,0.327)(0,1)
\or\qbezier(1,0)(0.179,0.311)(-0.5,0.866)
\or\qbezier(1,0)(0.0536,0.2)(-0.866,0.5)
\or\put(1, 0){\line(-2, 0){2}}
\or\qbezier(1,0)(0.0536,-0.2)(-0.866,-0.5)
\or\qbezier(1,0)(0.179,-0.311)(-0.5,-0.866)
\or\qbezier(1,0)(0.327,-0.327)(0,-1)
\or\qbezier(1,0)(0.45,-0.26)(0.5,-0.866)
\or\qbezier(1,0)(0.516,-0.138)(0.866,-0.5)
\fi
\or\ifcase#2\or
\or\qbezier(0.866,0.5)(0.378,0.378)(0.5,0.866)
\or\qbezier(0.866,0.5)(0.26,0.45)(0,1)
\or\qbezier(0.866,0.5)(0.12,0.446)(-0.5,0.866)
\or\qbezier(0.866,0.5)(0,0.359)(-0.866,0.5)
\or\qbezier(0.866,0.5)(-0.0536,0.2)(-1,0)
\or\put(0.866, 0.5){\line(-5, -3){1.73}}
\or\qbezier(0.866,0.5)(0.146,-0.146)(-0.5,-0.866)
\or\qbezier(0.866,0.5)(0.311,-0.179)(0,-1)
\or\qbezier(0.866,0.5)(0.446,-0.12)(0.5,-0.866)
\or\qbezier(0.866,0.5)(0.52,0)(0.866,-0.5)
\fi
\or\ifcase#2\or\or
\or\qbezier(0.5,0.866)(0.138,0.516)(0,1)
\or\qbezier(0.5,0.866)(0,0.52)(-0.5,0.866)
\or\qbezier(0.5,0.866)(-0.12,0.446)(-0.866,0.5)
\or\qbezier(0.5,0.866)(-0.179,0.311)(-1,0)
\or\qbezier(0.5,0.866)(-0.146,0.146)(-0.866,-0.5)
\or\put(0.5, 0.866){\line(-3, -5){1}}
\or\qbezier(0.5,0.866)(0.2,-0.0536)(0,-1)
\or\qbezier(0.5,0.866)(0.359,0)(0.5,-0.866)
\or\qbezier(0.5,0.866)(0.446,0.12)(0.866,-0.5)
\fi
\or\ifcase#2\or\or\or
\or\qbezier(0,1.)(-0.138,0.516)(-0.5,0.866)
\or\qbezier(0,1.)(-0.26,0.45)(-0.866,0.5)
\or\qbezier(0,1.)(-0.327,0.327)(-1,0)
\or\qbezier(0,1.)(-0.311,0.179)(-0.866,-0.5)
\or\qbezier(0,1.)(-0.2,0.0536)(-0.5,-0.866)
\or\put(0, 1){\line(0, -2){2}}
\or\qbezier(0,1.)(0.2,0.0536)(0.5,-0.866)
\or\qbezier(0,1.)(0.311,0.179)(0.866,-0.5)
\fi
\or\ifcase#2\or\or\or\or
\or\qbezier(-0.5,0.866)(-0.378,0.378)(-0.866,0.5)
\or\qbezier(-0.5,0.866)(-0.45,0.26)(-1,0)
\or\qbezier(-0.5,0.866)(-0.446,0.12)(-0.866,-0.5)
\or\qbezier(-0.5,0.866)(-0.359,0)(-0.5,-0.866)
\or\qbezier(-0.5,0.866)(-0.2,-0.0536)(0,-1)
\or\put(-0.5, 0.866){\line(3, -5){1}}
\or\qbezier(-0.5,0.866)(0.146,0.146)(0.866,-0.5)
\fi
\or\ifcase#2\or\or\or\or\or
\or\qbezier(-0.866,0.5)(-0.516,0.138)(-1,0)
\or\qbezier(-0.866,0.5)(-0.52,0)(-0.866,-0.5)
\or\qbezier(-0.866,0.5)(-0.446,-0.12)(-0.5,-0.866)
\or\qbezier(-0.866,0.5)(-0.311,-0.179)(0,-1)
\or\qbezier(-0.866,0.5)(-0.146,-0.146)(0.5,-0.866)
\or\put(-0.866, 0.5){\line(5, -3){1.73}}
\fi
\or\ifcase#2\or\or\or\or\or\or
\or\qbezier(-1,0)(-0.516,-0.138)(-0.866,-0.5)
\or\qbezier(-1,0)(-0.45,-0.26)(-0.5,-0.866)
\or\qbezier(-1,0)(-0.327,-0.327)(0,-1)
\or\qbezier(-1,0)(-0.179,-0.311)(0.5,-0.866)
\or\qbezier(-1,0)(-0.0536,-0.2)(0.866,-0.5)
\fi
\or\ifcase#2\or\or\or\or\or\or\or
\or\qbezier(-0.866,-0.5)(-0.378,-0.378)(-0.5,-0.866)
\or\qbezier(-0.866,-0.5)(-0.26,-0.45)(0,-1)
\or\qbezier(-0.866,-0.5)(-0.12,-0.446)(0.5,-0.866)
\or\qbezier(-0.866,-0.5)(0,-0.359)(0.866,-0.5)
\fi
\or\ifcase#2\or\or\or\or\or\or\or\or
\or\qbezier(-0.5,-0.866)(-0.138,-0.516)(0,-1)
\or\qbezier(-0.5,-0.866)(0,-0.52)(0.5,-0.866)
\or\qbezier(-0.5,-0.866)(0.12,-0.446)(0.866,-0.5)
\fi
\or\ifcase#2\or\or\or\or\or\or\or\or\or
\or\qbezier(0,-1.)(0.138,-0.516)(0.5,-0.866)
\or\qbezier(0,-1.)(0.26,-0.45)(0.866,-0.5)
\fi
\or\ifcase#2\or\or\or\or\or\or\or\or\or\or
\or\qbezier(0.5,-0.866)(0.378,-0.378)(0.866,-0.5)
\fi\fi\fi\fi}
\def\FullChord[#1,#2]{
\Endpoint[#1]
\Endpoint[#2]
\Arc[#1]
\Arc[#2]
\Chord[#1,#2]
}
\def\EndChord[#1,#2]{
\Endpoint[#1]
\Endpoint[#2]
\Chord[#1,#2]
}
\def\Picture#1{
\begin{picture}(2,1)(-1,-0.167)
#1
\end{picture}
}
\def\DottedChordDiagram[#1,#2]{
\Picture{\DottedCircle \FullChord[#1,#2]}
}
\def\u#1{ \underline{#1} }
\def\uu#1{ \underline{\underline{#1}} }

\def\ux{ {\u{x}} }
\def\uy{ {\u{y}} }

\def\ual{ {\u{\a}} }

\def\um{ {\u{m}} }

\def\uube{ {\uu{\beta}} }

\def\uuy{ {\uu{y}} }

\def\hcK{ \hat{\cK} }
\def\hcL{ \hat{\cL} }
\def\hcLp{ \hat{\cL}\p }

\def\chq{ \check{q} }
\def\chK{ \check{K} }
\def\chh{ \check{h} }
\def\cha{ \check{\a} }
\def\chb{ \check{\b} }
\def\chgamma{ \check{\gamma} }

\def\uchal{ \u{\cha} }
\def\uuchbe{ \uu{\chb} }

\def\fq{ q }
\def\fK{ K }
\def\fh{ h }

\def\RHS{R.H.S.}
\def\RHS{$\IQ$HS}

\def\prb#1{ \{ #1 \} }

\def\brlst#1#2#3{ \{ #1_{#2},\ldots,#1_{#3} \} }
\def\brclst#1#2#3{ ( #1_{#2},\ldots,#1_{#3} ) }

\def\szbK{ \sum_{0<\chb<\chK} }
\def\szubKo{ \sum_{0<\uuchbe<K\atop \uuchbe - {\rm odd} } }
\def\szbmK{ \sum_{-\chK<\chb\leq \chK\atop \chb - {\rm odd} } }
\def\szbKo{ \sum\limits_{0<\chb<\chK\atop \chb - {\rm odd} } }
\def\spzbKo{ {\sum\limits_{0<\chb<\chK\atop \chb - {\rm odd} } }
^{\!\!\!\prime} }
\def\sgmKh{ \sum_{-(\chK+1)/2<\chgamma\leq (\chK-1)/2} }
\def\sgoK{ \sum_{\chgamma=1}^{\chK} }
\def\snzKmt{ \sum_{n=0}^{K-2} }
\def\snzi{ \sum_{n=0}^{\infty} }
\def\snzN{ \sum_{n=0}^N }
\def\sumnz{ \sum_{\um,n\geq 0 \atop |\um|\leq n} }
\def\summnz{ \sum_{ {\um,m,n\geq 0\atop |\um|+m\leq n} \atop
        m\leq n/2} }
\def\szumn{ \sum_{\um\geq 0\atop |\um|\leq n} }
\def\sjoL{ \sum_{j=1}^L }

\def\smnz{ \sum_{m,n\geq 0} }
\def\smo{ \sum_{m\geq 1} }
\def\smz{ \sum_{m\geq 0} }
\def\snzmo{ \sum_{n=0}^{m+1} }
\def\sokLnj{ \sum_{1\leq k\leq L \atop k\neq j} }

\def\sjkoLp{ \sum_{j,k=1}^{L\p} }
\def\smnzN{ \sum_{m,n\geq 0\atop m+n\leq N} }

\def\pjoL{ \prod_{j=1}^{L} }

\def\Jbas#1#2{ J_{#1}(#2;\fq) }
\def\Jbasch#1#2{ J_{#1}(#2;\chq) }
\def\JuahL{ \Jbas{\uchal}{\hcL} }
\def\JuchahL{ \Jbasch{\uchal}{\hcL} }

\def\JuchaphL{ \Jbasch{\uchal\p}{\hcL} }
\def\JuaubhLhLp{ \Jbas{\ual,\uube}{\hLhLp} }
\def\JchuaubhLhLp{ \Jbasch{\uchal,\uuchbe}{\hLhLp} }

\def\Zbas#1#2#3{ Z^{#1}_{#2}(#3;K) }
\def\Zbasch#1#2#3{ Z^{#1}_{#2}(#3;\chK) }
\def\ZuaML{ \Zbas{}{\ual}{M,\cL} }
\def\ZuaMhL{ \Zbasch{}{\uchal}{M,\hcL} }
\def\ZcuaMhL{ \Zbasch{(c)}{\uchal}{M,\hcL} }

\def\Zchsbas#1{ Z(#1;\chK) }

\def\ZchuabMhLK{ \Zbasch{}{\uchal,\chb}{M,\hLK} }

\def\ZchuabMhLhK{ \Zbasch{}{\uchal,\chb}{M,\hLhK} }
\def\ZpM{ \Zbas{\prime}{}{M} }
\def\ZchpM{ \Zbasch{\prime}{}{M} }
\def\ZM{ Z(M;K) }
\def\ZchM{ Z(M;\chK) }
\def\ZcM{ Z\bc(M;K) }
\def\ZMt{ Z(M;3) }
\def\ZMtc{ \overline{\ZMt} }

\def\Ztrbas#1#2{ \Zbas{{\rm (tr)}}{#1}{#2} }
\def\Ztrbasch#1#2{ \Zbasch{{\rm (tr)}}{#1}{#2} }
\def\ZtruabMhLhK{ \Ztrbas{\ual,\b}{M,\hLhK} }
\def\ZtruchabMhLhK{ \Ztrbas{\uchal,\b}{M,\hLhK} }
\def\ZtruabMhLK{ \Ztrbas{\ual,\b}{M,\hLK} }
\def\ZtruaMhL{ \Ztrbas{\ual}{M,\hcL} }
\def\ZtruchaMhL{ \Ztrbasch{\uchal}{M,\hcL} }
\def\ZtruchaMhLm{ \Ztrbas{\uchal}{M,\hcL} }
\def\ZtruaML{ \Ztrbas{\ual}{M,\cL} }
\def\ZtruaMphL{ \Ztrbas{\ual}{M\p,\hcL} }

\def\ZtruchaMphL{ \Ztrbas{\uchal}{M\p,\hcL} }
\def\ZtrM{ \Ztrbas{}{M} }
\def\Zpbas#1#2{ \Zbas{\prime}{#1}{#2} }
\def\Zpbasch#1#2{ \Zbasch{\prime}{#1}{#2} }
\def\ZpuaMhL{ \Zpbas{\ual}{M,\hcL} }
\def\ZpuchaMhL{ \Zpbasch{\uchal}{M,\hcL} }

\def\ZpuchaMphL{ \Zpbasch{\uchal}{M\p,\hcL} }

\def\ZpuchabMhLhK{ \Zpbasch{\uchal,\chb}{M,\hLhK} }
\def\ZM{ Z(M;K) }
\def\ZtrLp{ \Ztrbas{}{\Lpo} }
\def\ZpLp{ \Zpbasch{}{\Lpo} }

\def\XhcLp#1{ X(\hcLp;#1) }
\def\XhcLpk{ \XhcLp{\kappa} }
\def\XcM{ X\bc(M;K) }

\def\Dbas#1#2{ D_{#1}(M,#2) }

\def\Pbas#1#2{ P_{#1}(M,#2) }
\def\Pebas#1{ P_{#1}(M) }
\def\PeMpbas#1{ P_{#1}(M\p) }
\def\PuanhL{ \Pbas{\ual\sepx n}{\hcL} }
\def\PuchanhL{ \Pbas{\uchal\sepx n}{\hcL} }
\def\Pwbas#1#2{ P^{\vee}_{#1}(M,#2) }

\def\PwuchanhL{ \Pwbas{\uchal\sepx n}{\hcL} }

\def\PuchanMhLKv#1{ \Pbas{\uchal\sepx n}{\hcL,\cK;#1} }

\def\PuchanMhLKt{ \PuchanMhLKv{t} }

\def\PuchanMhLKqb{ \PuchanMhLKv{q^\b} }

\def\dMKbas#1{ d_{#1}(M,\cK) }
\def\dzz{ \dMKbas{0;0} }
\def\doz{ \dMKbas{1;0} }
\def\dzo{ \dMKbas{0;1} }
\def\dmnK{ \dMKbas{m;n} }

\def\DA{ \Delta_{\rm A} }
\def\DAtnMKqb{ \DA^{2n+1}(M,\cK;q^\b) }
\def\DAMKv#1{ \DA(M,\cK;#1) }
\def\DAMKt{ \DAMKv{t} }
\def\DAddMK{ \DA^{\prime\prime}(M,\cK) }

\def\dbas#1#2{ d_{#1}(M,#2) }
\def\dumnhL{ \dbas{\um\sepx n}{\hcL} }

\def\duchamnMhLK{ \dbas{\uchal,m\sepx n}{\hcL,\cK} }

\def\dchv{ d^{\vee}_{\uchal,m\sepx n}(M,\hcL,\cK) }

\def\ZZqhpm{ \ZZ[q^{1/2},q^{-1/2}] }
\def\ZZtpm{ \ZZ[t,t^{-1}] }
\def\ZZqpm{ \ZZ[q,q^{-1}] }
\def\ZZcq{ \ZZ[\chq] }

\def\ZZhordMi{ \ZZ[1/2,1/\ordhM] }
\def\ZZhordMpi{ \ZZ[1/2,1/\ordhozbas{M\p}] }
\def\ZZordMih{ \ZZ[1/\ordhM]\,[[h]] }
\def\ZZtohord{ \ZZ[t^{1/o},t^{-1/o},1/2,1/\ordhM] }
\def\ZZqhord{ \ZZ[q,q^{-1},1/2,1/\ordhM] }
\def\ZZhpi{ \ZZ[1/p]\,[[h]] }

\def\ZZK{ \ZZ_{\chK} }
\def\ZZKcq{ \ZZK[\chq] }

\def\QchK{ \ZZ_{(\chK)} }

\def\bbrif#1#2#3#4#5{ \left\{ \begin{array}{cl}
        #1 &\mbox{if $#2$}\\
        #3 &\mbox{if $#4$}{#5}
        \end{array} \right. }

\def\eqmod#1#2#3{ #1\equiv #2\,\,\mod{#3} }

\def\hozbas#1{ H_1(#1,\ZZ) }
\def\homz{ \hozbas{M} }
\def\ordhozbas#1{ |\hozbas{#1}| }
\def\ordhM{ \ordhozbas{M} }
\def\ordhMh{ \ordhM^{1/2} }
\def\ordhMph{ \ordhMp^{1/2} }
\def\ordhMmh{ \ordhM^{-1/2} }
\def\ordhMp{ \ordhozbas{M\p} }

\def\emtfs#1{ e^{-(3\pi i/4)\prosign(#1)} }
\def\emtfsLp{ \emtfs{\hcLp} }
\def\emtfsp{ \emtfs{p} }

\def\emfs#1#2{ e^{(i\pi/4)(\kappa #2 1)\prosign(#1)} }
\def\emfsLp{ \emfs{\hcLp}{+} }
\def\emfsp{ \emfs{p}{-} }
\def\empfsp{ e^{-(i\pi/4)\prosign(p)} }

\def\qtfs#1{ q^{(3/4)\prosign(#1)} }
\def\qchtfs#1{ \chq^{(3/4)\prosign(#1)} }

\def\qchtfsLp{ \qchtfs{\hcLp} }
\def\qtfsp{ \qtfs{p} }

\def\qtfss#1{ \chq^{3\, 4\s\prosign(#1)} }
\def\qtfssLp{ \qtfss{\hcLp} }
\def\qtfssp{ \qtfss{p} }

\def\mos#1{ (-1)^{\prosign(#1)} }
\def\mosLp{ \mos{\hcLp} }

\def\tpi{ 2\pi i }

\def\qpmh#1{ q^{#1/2} - q^{-#1/2} }
\def\qchpmh#1{ \chq^{#1/2} - \chq^{-#1/2} }
\def\epmh#1{ e^{(i\pi/K)#1} - e^{-(i\pi/K)#1} }

\def\qmbtm{ (\qpmh{\b})^{2m+2} }
\def\chqmbtm{ (\qchpmh{\chb})^{2m+2} }

\def\chqmbto{ (\qchpmh{\chb})^{2m+1} }

\def\imLp{ i^{-L\p} }

\def\tchKmLph{ (2\chK)^{-L\p/2} }
\def\KmLph{ K^{-L\p/2} }
\def\chKmLph{ \chK^{-L\p/2} }
\def\eKmoph{ (8K)^{-1/2} }
\def\hLhLp{ \hcL\cup\hcLp }
\def\hLhK{ \hcL\cup\hcK }
\def\hLK{ \hcL\cup\cK }

\def\hLphK{ \hcLp\cup\hcK }

\def\anMchK{ a_n(M;\chK) }
\def\lnM{ \l_n(M) }

\def\wdg#1{ #1^{\vee} }

\def\lw{ \wdg{\l} }

\def\Legsbas#1{ \lrbc{#1\over K}_{\rm\!\! L} }
\def\Legsbasch#1{ \lrbc{#1\over \chK}_{\rm\!\! L} }
\def\LegsM{ \Legsbas{\ordhM} }
\def\LegschM{ \Legsbasch{\ordhM} }

\def\Legschp{ \Legsbasch{|p|} }

\def\LegschMp{ \Legsbasch{\ordhMp} }

\def\lcw{ \l_{\rm CW} }
\def\lcwM{ \lcw(M) }

\def\Keq{ \mathop{=}\limits_{ \chK-{\rm adic} } }
\def\Keq{ \mathop{=} }

\def\Ki{ K^{-1} }
\def\chKi{ \chK^{-1} }

\def\chKh{ \chK^{1/2} }
\def\Kmh{ K^{-1/2} }
\def\chKmh{ \chK^{-1/2} }

\def\uaumo{ \ual^{2\um+1} }
\def\bmo{ \b^{2m+1} }

\def\istbz{ \int_{[\b=0]} }
\def\intmzz{ \int_{-\infty}^{+\infty} }

\def\qrtrv{ {1\over 4} }
\def\hlfv{ {1\over 2} }
\def\Scsc{ S_{\rm CS}^{(c)} }
\def\XcuaMhL{ X^{(c)}_{\uchal}(M,\hcL;\chK) }

\def\s{ ^{\ast} }
\def\wve{ ^{\vee} }
\def\sip{ \prosign(p) }

\def\thpm{ t^{1/2} - t^{-1/2} }

\def\qfspb{ \chq^{4\s p\chb^2} }
\def\qfpb{ q^{(1/4)p\b^2} }
\def\qmmps{ \chq^{-m^2 p\s} }
\def\qmmp{ q^{-m^2/p} }
\def\qmbp{ (q^{m\b} + q^{-m\b}) }
\def\chqmbp{ (\chq^{m\chb} + \chq^{-m\chb}) }

\def\pah{ |p|^{1/2} }

\def\hi{ h^{-1} }
\def\chhi{ \chh^{-1} }

\def\xon{ x_1 }
\def\xtw{ x_2 }
\def\xth{ x_3 }
\def\xfo{ x_4 }

\def\Lpo{ L_{p,1} }

\def\bc{ ^{(c)} }

\def\skzi{ \sum_{k=0}^\infty }
\def\fmn{ {m \over n} }

\def\sepx{ ; }

\def\nssN{ \remkN{(1/n)\wve} }

\begin{document}

\begin{titlepage}
\centerline{\hfill                 math.QA/9806075}
\vfill
\begin{center}

{\large \bf
On $p$-adic properties of the Witten-Reshetikhin-Turaev invariant
} \\

\bigskip
\centerline{L. Rozansky\footnote{
This work was supported by NSF Grant DMS-9704893}
}

\centerline{\em Department of Mathematics, University of Illinois}
\centerline{\em Chicago, IL 60608, U.S.A.}
\centerline{{\em E-mail address: rozansky@math.uic.edu}}

\vfill
{\bf Abstract}

\end{center}
\begin{quotation}

We prove the Lawrence conjecture about $p$-adic convergence of the
series of Ohtsuki invariants of a rational homology sphere to its
$SO(3)$ Witten-Resheti-\linebreak khin-Turaev invariant.
Our proof is based on the
surgery formula for Ohtsuki series and on the properties of the
expansion of the colored Jones polynomial of a knot in powers of
$q-1$ and $q^\a-1$, $\a$ being the color of the knot.

\end{quotation} \vfill \end{titlepage}

\pagebreak

\nsection{Introduction}
\label{s1}
\hyphenation{Re-she-ti-khin}
\hyphenation{Tu-ra-ev}

A quantum invariant of 3-manifolds (WRT invariant) discovered by
E.~Witten\cx{Wi} and by N.~Reshetikhin, V.~Turaev\cx{ReTu}, is an
extension of the Jones polynomial from links in $S^3$ to links in any
compact oriented 3-manifold $M$. For a link
$\cL\subset M$ with $L$ components,
 the WRT invariant $\ZuaML$ is a complex number
depending on positive integer parameters $K$ and
$\ual = \brclst{\a}{1}{L}$ such that $1\leq \a_j\leq K-1$. The WRT
invariant is quite effective in distinguishing between the
topologically inequivalent links and 3-manifolds. However its
relation to the `classical' topology remains mysterious, especially
if we use the mathematically rigorous definition of\cx{ReTu}, which
relies on the representation theory of quantum group $su_q(2)$
rather than the non-rigorous path integral definition
of\cx{Wi}.


One could try to determine the
topological content of the WRT invariant $\ZM$
independently for every value of $K$. R.~Kirby
and P.~Melvin have tried this approach in\cx{KiMe}, but they
succeeded in relating $\ZM$ to classical invariants only for few
small values of $K$. An alternative approach is to study how $\ZM$
depends on $K$ in a hope of identifying the topological
invariants which parametrize this dependence. Two different methods of
achieving this goal were developed -- a number-theoretic and an
analytic one. We showed in\cx{Ro5} that
if $M$ is a rational homology sphere (\RHS) then, surprisingly, both
methods produce the same sequence of invariants.

The number-theoretic method was advanced by H.~Murakami and
T.~Ohtsuki.
They worked with the modified WRT invariant
$\ZpM$ introduced by Kirby and Melvin\cx{KiMe}.
%
Murakami proved\cx{Mu1},\cx{Mu2} (see also\cx{MR}) that
for a \RHS\ $M$ and for an odd prime $K$,
$\ZpM$ belongs to the cyclotomic ring $\ZZ[q]$, $q=\exp(2\pi i/K)$.
Therefore one can write
\qq
\ZpM = \snzKmt a_n(M;K)\,(q-1)^n,\qquad
\mbox{where $a_n(M;K)\in\ZZ$}.
\label{v5}
\qqq
%
The coefficients $a_n(M;K)$ still depend on $K$, however Ohtsuki
observed\cx{Oh1},\cx{Oh2}
that their remainders modulo $K$ are almost $K$-independent.
He showed that there exist $K$-independent invariants
$\lnM\in\IQ$, $n\geq 0$ such that
\qq
\eqmod{a_n(M;K)}{\LegsM\lnM}{K}\;\;\;\mbox{for $n\leq {K-3\over 2}$}.
\label{v6}
\qqq
Here $\ordhM$ is the order of the first homology of $M$ and
$\LegsM$ is the Legendre symbol.
If $\lnM$ is not integer,
then it should be viewed in \ex{v6} as a $K$-adic integer.
Murakami showed that the first two invariants $\l_0$ and $\l_1$ are
classical
%
\qq
\l_0(M) = 1/\ordhM,\qquad \l_1(M) = 3\lcwM/\ordhM,
\label{v09}
\qqq
where $\lcwM$ is the Casson-Walker invariant of $M$.

R.~Lawrence\cx{La} conjectured that $\lnM$ have much more
control over $\ZpM$, than that
displayed by \ex{v6}. She suggested that
\qq
\lnM \in \ZZhordMi
\label{v8}
\qqq
and that
the series
$
\LegsM\snzi \lnM\,(q-1)^n
$
converges $K$-adically to $\ZpM$
in the $K$-adic completion $\ZZ_K[q]$ of the cyclotomic ring
$\ZZ[q]$,
that is,
\qq
\ZpM = \LegsM\snzi \lnM\,(q-1)^n.
\label{v9}
\qqq
Lawrence verified her conjecture for Seifert manifolds constructed by
Dehn's surgery on $(2,m)$ torus knots.

%
%
%

The goal of this paper is to prove the Lawrence conjecture for a
general \RHS\ Our proof is based on a surgery formula for the
invariants $\lnM$. This formula has first appeared\cx{Ro1} in the
study of the asymptotic properties of $\ZM$ at large values of $K$.
The same formula is relevant for the number-theoretic properties of
$\ZpM$ because
the sum in the \rhs of \ex{v9} considered as a
formal power series in $(q-1)$, represents the trivial connection
contribution to $\ZM$.

At the `physical' level of rigor,
the asymptotic properties of $\ZM$ follow from
Witten's formula\cx{Wi} which presents this quantum invariant
as a path integral over
$SU(2)$ connections in the (trivial) $SU(2)$ bundle over $M$.
%
At large values
of $K$ the integral can be evaluated in the stationary phase
approximation, the stationary phase points being the flat $SU(2)$
connections. Therefore Witten conjectured that $\ZM$ splits into
a sum of contributions coming from connected components
$\cM\bc$ of the moduli space $\cM$ of flat $SU(2)$ connections on $M$
\qq
\ZM \simeq \sum_c \ZcM,
\label{v1}
\qqq
where $\simeq$ denotes asymptotic convergence.
Each contribution $\ZcM$ is proportional to a `classical
exponential'
\qq
\ZcM = e^{(iK/2\pi)\Scsc} \XcM.
\label{v2}
\qqq
Here $\Scsc$ is the Chern-Simons invariant of connections
in the component $\cM\bc\subset \cM$
while $\XcM$ is an asymptotic series in
powers of $K^{-1}$ whose coefficients are invariants of $M$. These
coefficients are independent of $K$ and therefore they can be expected
to be related to classical topology.
If $M$ is a \RHS, then the trivial $SU(2)$
connection is an isolated point in the moduli space $\cM$. Therefore
according to path integral predictions, it should yield a distinct
contribution to $\ZM$ of the form
\qq
\ZtrM = \ordhMmh \snzi D_n(M)\,K^{-n},
\label{v3}
\qqq
where $D_n(M)$
are complex-valued invariants of $M$.

Since the asymptotic properties\rx{v1},\rx{v2} of the WRT invariant
have not been proven yet, we had to
choose a different way to work with
$\ZtrM$. Following the approach of Reshetikhin and Turaev\cx{ReTu}
to the WRT invariant, we
{\em defined\/} the invariant $\ZtrM$ of a \RHS\ $M$
constructed by a surgery on a link $\cL\subset S^3$
as a formal power series\rx{v3} whose
coefficients $D_n(M)$ are expressed in terms of the coefficients of
the colored Jones polynomial of $\cL$ through special surgery
formulas\cx{Ro1}-\cite{Ro6}. We call $\ZtrM$ defined by these
formulas {\em the TCC invariant}, TCC being the abbreviation of
the trivial connection contribution.
We proved that $\ZtrM$ is well-defined by our
 surgery formulas, that is, it does not depend on
the choice of a surgery link used
to construct $M$. Then we gave a path
integral explanation of why
our invariant should indeed represent the trivial connection
contribution to the path integral of\cx{Wi}.
We also proved\cx{RoS1},\cx{RoS2}
that the WRT invariants of Seifert homology spheres
have the
asymptotic properties\rx{v1},\rx{v2} and that the surgery formula for
the TCC invariant of these manifolds yields the trivial
connection contribution to their WRT invariants in
the sense of \eex{v1},\rx{v2}.

Our surgery formulas for $\ZtrM$
showed that, as we hoped, the first coefficients $D_n(M)$ in the sum
of \ex{v3}
are
classical topological invariants of $M$
\qq
D_0(M) = 1/\ordhM,\qquad D_1(M) = 6\pi i\lcwM/\ordhM.
\label{v4}
\qqq
%
There is an apparent similarity between \eex{v4} and\rx{v09}.
In\cx{Ro5} we showed that for $n\geq 2$ the invariants $D_n(M)$ and
$\lnM$ are also
essentially equivalent. More precisely, we proved the equation
\qq
\ZtrM = \snzi \lnM\,(e^{2\pi i/K} - 1)^n,
\label{v10}
\qqq
which should be understood as a relation between two
formal power series in $K^{-1}$.
In view of relation\rx{v10}, the Lawrence conjecture\rx{v9} has
a deeper meaning: it means that although the
TCC invariant
$\ZtrM$ is only a part of the
total WRT invariant $\ZM$, still it determines $\ZM$ for all prime
$K$.

Our proof of \ex{v10} as well as an alternative proof of the Murakami
and Ohtsuki results
\rx{v5}
and\rx{v09} which we provided
in\cx{Ro5}, was based on the propreties of the
expansion of the colored Jones polynomial in
powers of $(q-1)$ and colors
(this expansion was first considered
in\cx{MeMo}). We used a previously established bound\cx{Ro4} on the
powers of colors versus the powers of $(q-1)$ which exists for
algebraically split links, \ie links whose linking numbers are zero
(recently T.~Le\cx{Le} has used the same idea in order to prove the
analogs of \eex{v5} and\rx{v09} for $SU(N)$ WRT invariants). However,
the bound itself did not allow us to relate $\lnM$ to $\ZpM$ beyond
\ex{v09}.

In\cx{Ro9} we proved that the expansion of
the colored Jones polynomial of a knot in $S^3$ in powers of $(q-1)$
and $\a$ ($\a$ being the color of the knot)
can be rewritten as
an expansion in powers of $(q-1)$ and $(q^\a - 1)$ ($\a$ being the
color of the knot) with integer coefficients. The integrality of
those coefficients allowed us\cx{Ro8} to give a simple proof of
Lawrence's conjecture for a \RHS, constructed by a surgery on a knot
in $S^3$. In\cx{Ro11} we generalized the results of\cx{Ro9} to a knot
in a \RHS\ In this paper we will use the method of\cx{Ro8} in order
to prove the Lawrence conjecture in the general case.


\nsection{Notations, background and statement of results}
\label{s2*}

First of all, let us introduce our multi-index, number-theoretic and
topological notations.

We use multi-index notations for the colors of
links.
For a link $\cL$ with $L$ components we denote
%
$\ux = \brclst{x}{1}{L}$
%
and
\qq
\begin{array}{c@{\qquad} c@{\qquad} c}
\ux+\uy = (x_1+y_1,\ldots,x_L+y_L), &
y\ux  =  \brclst{yx}{1}{L}, &
\ux\uy  =  (x_1 y_1,\ldots,x_L y_L),\\[2ex]
x^{\uy}  =  (x^{y_1},\ldots,x^{y_L}), &
\ux^y  =  \brclst{x^y}{1}{L}, &
\ux^{\uy} = \pjoL x_j^{y_j},\\[2ex]
\prb{f(\ux)}  =  \pjoL f(x_j), &
|\ux|  =  \sjoL x_j. &
\end{array}
\label{va.2}
\qqq
When we consider two links $\cL$ and $\cL\p$ simultaneously,
we also use a multi-index notation
\qq
\uuy = \brclst{y}{1}{L\p}
\qqq
for the colors of the $L\p$-component link $\cL\p$.
We denote concatenation of multi-indices as
\qq
\ux,x = (x_1,\ldots, x_L,x),\qquad
\ux,\uuy = (x_1,\ldots, x_L, y_1,\ldots, y_{L\p}).
\qqq
Also
%
$\ux = x$ means $x_j=x$ for all $1\leq j\leq L$.

We will use three variables $\fK$, $\fq$ and $\fh$ which are
related by
\qq
\fq = e^{2\pi i/\fK},\qquad \fh = \fq -1.
\label{v2.1}
\qqq
A rational power of $q$ is defined as
\qq
q^r = e^{2\pi ir/K},\quad r\in\IQ.
\qqq
In many cases we will have to specialize $\fK$ to a positive
integer. In order to indicate this explicitly we use `checked'
symbols
\qq
\chK \in\ZZ_+,\qquad \chq = e^{2\pi i/\chK},\qquad
\chh = \chq - 1.
\label{v2.2}
\qqq
Also while $\a,\b$ are variables, $\cha$ and $\chb$ stand
for positive integers.
These elaborate notations will
help us to avoid
confusion between the formal parameters and their
$K$-adic and cyclotomic versions.

\def\remk#1#2{ \left[#1\right]_{#2} }
\def\remkN#1{ \remk{#1}{N} }

For $\chK$ being prime, $\ZZK$ denotes the ring of $\chK$-adic
integers, $\ZZcq$ is the cyclotomic ring and $\ZZKcq$
is its $K$-adic completion. Note that since we set
$\chq = e^{2\pi i/\chK}$, then $\ZZcq$ for us is a particular subring
of $\IC$.

Denote by $\QchK\subset \IQ$ the ring of
rational numbers whose denominators are not divisible by $\chK$.
There is a natural embedding
\qq
{}\wve:\,\QchK\rightarrow \ZZK,\qquad \l \mapsto \lw,
\label{map1}
\qqq
which maps fractions into power series in $\chK$.
Whenever we use a notation $\lw$, we assume that
$\l\in\QchK$.
The homomorphism \rx{map1} can be extended from $\QchK$ to formal
power series with coefficients in $\QchK$
\qq
{}\wve:\,\QchK[[h]]\rightarrow\ZZKcq
\label{map2}
\qqq
by setting $h\mapsto \chh$. Indeed, since
\qq
(\chq^{\chK} - 1)/(\chq -1) = 0,
\label{v2.3}
\qqq
then
\qq
\chh^{\chK-1} = \chK x,\qquad \mbox{for some $x\in\ZZcq$.}
\label{v2.4}
\qqq
As a result,
any series of the form $\snzi a_n\,\chh^n$, $a_n\in\ZZK$
has a $\chK$-adic limit in $\ZZKcq$
\qq
\snzi a_n\,\chh^n \Keq A\in\ZZKcq,
\label{vb2.1}
\qqq
which means by definition that for any $N_0>0$ there exists
$N_0\p$ such that for any $N>N_0\p$
\qq
A = \snzN a_n\,\chh^n + \chK^{N_0+1}x,
\qquad x\in \ZZKcq.
\label{1.18}
\qqq
Therefore the image of a formal power series
$S(h)=\snzi a_n\,h^n$, $a_n\in\QchK$ under the homomorphism\rx{map2}
can be defined as the $K$-adic limit of the corresponding series
$\snzi a_n\wve\,\chh^n$
\qq
{}\wve:\, S(h) \mapsto S\wve \Keq \snzi a_n\wve\chh^n \in\ZZKcq.
\label{map3}
\qqq

We will need a few more number-theoretic notations.
$\Legsbasch{x}$ denotes a Legendre symbol of $x$. If $x$ is an
integer not divisible by $\chK$, then
$\Legsbasch{x}=1$ if there exists $y\in\ZZ$ such that
$\eqmod{x}{y^2}{\chK}$, and $\Legsbasch{x}=-1$ otherwise.
For a prime
$\chK$ and an integer $n$ not divisible by a prime $\chK$,
let $n\s$ denote any integer such that
\qq
\eqmod{nn\s}{1}{\chK}.
\label{1.37}
\qqq
Also denote for an odd $\chK$
\qq
\kappa =
\bbrif{ 1 }{\eqmod{\chK}{1}{4}}
{ -1 }{\eqmod{\chK}{-1}{4}}. 
\label{1.38}
\qqq
Finally, for $x\in\ZZK$, we denote by
$\remkN{x}\in\ZZ$ the remainder of $x$ modulo $\chK^{N+1}$
\qq
\eqmod{\remkN{x}}{x}{\chK^{N+1}}.
\label{1.12}
\qqq

Now let us fix some standard topological notations.
Let $\cK$ be a knot in a
3-manifold $M$. A
{\em meridian\/} of $\cK$ is a simple cycle on the
boundary of its tubular neighborhood, which is contractible through
that neighborhood. A {\em parallel\/} is a cycle on the same boundary
which has a unit intersection number with the meridian. A knot is
called {\em framed\/}
if we made a choice of
its parallel. A link is framed if all of its components are framed.
We denote framed knots and links by putting a hat on top of their
symbols: $\hcK, \hcL$.

Suppose that a framed knot $\hcK$ is of finite order as
an element of $\homz$ (we will denote the order of $\hcK$ as $o$).
Then its {\em self-linking number\/} $p$ is the linking number
between the knot and its parallel,
$p\in\ZZ/o$.
{\em Dehn's
surgery\/} on a framed knot $\hcK\subset M$ is a transformation of
cutting out a tubular neighborhood of $\hcK$ and gluing it back in
such a way that the meridian of the tubular neighborhood matches the
parallel left on the boundary of the knot complement. A result of
this procedure is a new 3-manifold $M\p$.

In order to shorten our formulas, we will denote the order of the
first homology of a \RHS\ $M$ as $h_1(M)$:
\qq
h_1(M) = \ordhM.
\qqq

\def\ordhozbas#1{ h_1(#1) }
\def\ordhM{ \ordhozbas{M} }
\def\ordhMh{ h_1^{1/2}(M) }
\def\ordhMph{ h_1^{1/2}(M\p) }
\def\ordhMmh{ h_1^{-1/2}(M) }
\def\ordhMpmh{ h_1^{-1/2}(M\p) }
\def\ordhMp{ \ordhozbas{M\p} }
\def\ordhMmt{ h_1^{-2}(M) }



Let us recall the quantum invariants appearing in
3-dimensional topology.
The colored Jones polynomial
\qq
&\JuahL \in \ZZqhpm,
\label{1.1}\\
&\JuahL \in \ZZqpm\;\;\;\mbox{if $\uchal$ are odd}
\label{1.2}
\qqq
is an invariant of a framed $L$-component link $\hcL\subset S^3$.
The positive integers $\uchal$ are called `colors', they are assigned
to components of $\hcL$.
The
WRT invariant is an invariant of a framed link
$\hcL$ in a 3-manifold $M$. It depends on a positive
integer $\chK$, and we denote it as $\ZuaMhL$. N.~Reshetikhin and
V.~Turaev defined the WRT invariant by a surgery formula and verified
that $\ZuaMhL$ does not depend on the choice of a surgery link.
\begin{definition}[N.~Reshetikhin, V.~Turaev\cx{ReTu}]
\label{t1.1}
If $M=S^3$, then
\qq
&\Zbasch{}{\uchal}{S^3,\hcL} = \JuchahL,
\label{1.3}
\qqq
If $M$ is constructed by Dehn's surgery on a framed $L\p$-component
link $\hcLp\subset S^3$, then
\qq
\ZuaMhL & = & \imLp\,\tchKmLph\,\emtfsLp\,\qchtfsLp
\label{1.5}\\
&&\qquad\times\szbK \prb{\qchpmh{\uuchbe}}\,
\JchuaubhLhLp,
\nonumber
\qqq
where $\prosign(\hcLp)$ is
the signature of the linking matrix of $\hcLp$.
\end{definition}
For an empty link $\cL$ we denote $\ZuaMhL$ simply as $\ZchM$.
Note that we use the normalization in which
\qq
\Zchsbas{S^3} = 1.
\label{sp.1}
\qqq

The WRT invariant (and all other quantum invariants of links that we
will consider in this paper) has a simple dependence on the choice of
framing of $\hcL$. As a result, if a component $\hcK$ of a link
$\hLhK$ has a finite order as an element of $\homz$ so that it
has a well-defined self-linking number $p$, then we can
introduce an invariant which does not depend on the framing of $\hcK$
\qq
\ZchuabMhLK = q^{-(1/4)p(\chb^2-1)}\,\ZchuabMhLhK.
\label{1.6}
\qqq
If $\hcK$ is homologically trivial, then $\cK$ can be interpreted as
a knot with zero self-linking number, otherwise $\cK$ is just a
symbolic notation for the framing-independent invariant $\ZchuabMhLK$.

For an odd number $\chK$, R.~Kirby and P.~Melvin introduced the
$SO(3)$ WRT invariant $\ZchpM$ which has a simple relation to $\ZchM$
\qq
\begin{array}{cl}
\ZM = \ZMt\,\ZpM &\mbox{if $\eqmod{K}{-1}{4}$},\\
\ZM = \ZMtc\,\ZpM &\mbox{if $\eqmod{K}{1}{4}$}.
\end{array}
\label{1.7}
\qqq
%
H.~Murakami\cx{Mu1},\cite{Mu2}
and T.~Ohtsuki\cx{Oh1},\cite{Oh2}
proved that $\ZchpM$ satisfies the
following properties.

\begin{theorem}
\label{tMO}
If $\chK$ is an odd prime and $M$ is a \RHS\
such that the order of the first homology $\ordhM$ is not
divisible by $\chK$, then
\begin{itemize}
\item[(1)] [Murakami]
\qq
\ZchpM \in \ZZcq,
\label{1.8}
\qqq
%
so that
\qq
&\ZchpM = \snzKmt \anMchK\,\chh^n,\qquad \mbox{where $\anMchK\in\ZZ$;
}
\label{1.10}
\qqq
\item[(2)] [Ohtsuki]
There exists a sequence
of invariants $\lnM\in\IQ$, $n\geq 0$ such that
if $ n \leq {\chK - 3\over 2}$ then $\lnM\in\QchK$ and
\qq
\eqmod{\anMchK}{ \LegschM\,\lw_n(M)}{\chK},
\label{1.13}
\qqq
where $\Legsbasch{\cdot}$ is the Legendre symbol;

\item[(3)] [Murakami]
\qq
\l_0(M) = 1/\ordhM,\qquad \l_1(M) = 3\lcwM/\ordhM,
\label{1.14}
\qqq
where $\lcwM$ is the Casson-Walker invariant of $M$.

\end{itemize}
\end{theorem}

%
%
The goal of this paper is to prove the following conjecture made by
R.~Lawrence.
\begin{conjecture}[R.~Lawrence\cx{La}]
\label{t1.5}
For a \RHS\ M,
\qq
\lnM \in
\bbrif{\ZZ}{\ordhM = 1}{\ZZhordMi}{\ordhM>1}{}
\label{1.19}
\qqq
and if $\chK$ is an odd prime which does not divide $\ordhM$, then
\qq
\ZchpM \Keq \LegschM \snzi \lw_n(M)\,\chh^n
\quad\mbox{in $\ZZKcq$}.
\label{1.20}
\qqq
\end{conjecture}

The proof of the Lawrence conjecture will be based on the properties
of the invariant of \RHS\ which we call {\em the trivial connection
contribution\/} to the WRT invariant, or simply the TCC invariant. We
introduced it and studied its properties in\cx{Ro1}-\cite{Ro6}. The
TCC invariant of a \RHS\ $M$
is a formal power series in powers of $\Ki$
%
\qq
\ZtruaMhL = \ordhMmh \sumnz \Dbas{\um\sepx n}{\hcL}\,\uaumo\,K^{-n},
\label{1.21}
\qqq
%
whose coefficients
$\Dbas{\um,n}{\hcL}$ are
invariants of $M$ and $\hcL$.
The invariant\rx{1.21} can be defined by various equivalent surgery
formulas. Here we will use the definition which is provided by the
following
\begin{theorem}[\cite{Ro6}]
\label{t1.6}
There exists a unique invariant $\ZtruaMhL$ of
a framed link $\hcL$ in a \RHS\ $M$
(called the TCC invariant) which has
the form\rx{1.21} and satisfies the following three properties:
\begin{itemize}
\item[(1)]
If $M=S^3$, then \ex{1.21} coincides with the Melvin-Morton
expansion\cx{MeMo} of the colored Jones polynomial in powers of
$\Ki$ at fixed values of colors $\uchal$, that is
\qq
\Ztrbas{\uchal}{S^3,\hcL} = \JuchahL.
\label{1.22}
\qqq
\item[(2)]
Let $\hcL$ be a framed $L$-component link and $\hcK$ be a framed knot
with self-linking number $p$ in a \RHS\ $M$. Then the expression
%
%
\qq
&\ZtruabMhLK =
q^{-p(\b^2-1)/4}\,\ZtruabMhLhK,
\label{1.23}
\qqq
does not depend on the choice of the framing of $\hcK$
(that is, $\ZtruabMhLK$ is well-defined for unframed knots).
Moreover, if
$m>n/2$, then $\Dbas{\um,m,n}{\hLK}=0$ or, in other words,
%
\qq
&\ZtruabMhLK = \ordhMmh
\hspace{-5pt}\summnz\hspace{-5pt}
\Dbas{\um,m\sepx n}{\hLK}\,\uaumo\,\bmo\,
K^{-n}.
\label{1.24}
\qqq
\item[(3)]
If a \RHS\ $M\p$ is constructed by Dehn's surgery on a knot $\hcK$
with self-linking number $p$ in a \RHS\ $M$, then
\qq
\ZtruaMphL & = & - i\,\eKmoph\,\emtfsp\,\qtfsp
\label{1.25}\\
&&\qquad\times
\istbz  \ZtruabMhLhK\,(\qpmh{\b})\,d\b
\nonumber
\qqq
(here $\istbz$ denotes a stationary phase contribution of the point
$\b=0$), or more precisely,
%
%
%
\qq
\ZtruaMphL
& = & - i\,\eKmoph\,\emtfsp\,\qtfsp \ordhMmh
\label{1.26}\\
&&\qquad\times
\summnz
\Bigg(\Dbas{\um,m\sepx n}{\hLK}\,\uaumo\,K^{-n}
\nonumber\\
&&
\qquad\qquad\times
\left.\intmzz e^{(i\pi/2K) p (\b^2-1)}\,\bmo\,(\epmh{\b})\,d\b
\right).
\nonumber
\qqq

\end{itemize}
\end{theorem}
Equation\rx{1.26} is a well-defined relation between formal
power series. Indeed, since
\qq
K^{-1/2}
\intmzz e^{(i\pi/2K) p (\b^2-1)}\,\bmo\,(\epmh{\b})\,d\b
= O(K^{m})\quad\mbox{as $K\rightarrow \infty$,}
\label{sp.2}
\qqq
then only the coefficients $\Dbas{\um,m,n}{\hLK}$ with $n-m\leq n\p$
participate in the expression for $D_{\um,n\p}(M\p,\hcL)$. In view of
the bound $m\leq n/2$ in the sum of \ex{1.24}, this means that the
number of such coefficients $\Dbas{\um,m,n}{\hLK}$ is finite for any
fixed $n\p$.

We proved Theorem\rw{t1.6} in\cx{Ro6} by showing that
\ex{1.24} (which indicates that the removal of the self-linking
factor $q^{p(\b^2-1)/4}$
reduces the power of $\b$ versus the power of $\Ki$ in the expansion)
is consistent with
the surgery formula\rx{1.25} and
that the invariant $\ZtruaMhL$ does
not depend on the choice of the sequence of knot
 surgeries which leads from $S^3$ to
$M$. We also proved in\cx{Ro11} (Theorem~1.15 and Remark~1.17) that
if we switch in \ex{1.21} from powers of $\Ki$ to powers of
$h=e^{\tpi/K} - 1$ by writing
\qq
&\ZtruaMhL = \ordhMmh\snzi \PuanhL\,h^n,
\label{1.27}\\
&\PuanhL = \szumn \dumnhL\,\uaumo,
\label{1.28}
\qqq
then for fixed odd $\uchal$
\qq
\PuchanhL \in
\bbrif{\ZZ}{\ordhM = 1}{\ZZhordMi}{\ordhM>1}{}
\label{1.29}
\qqq
\begin{remark}
\label{t1.7}
\rm
In fact, we proved relation\rx{1.29} in\cx{Ro11} for a 0-framed link
$\cL$. However, since
\qq
\ZtruaMhL= q^{\qrtrv \sjoL l_{jj} (\a_j^2-1)}
\ZtruaML,
\qqq
$l_{jj}$ being self-linking numbers of the link components
$\hcL_j$, and since for odd $\uchal$
\qq
q^{\qrtrv \sjoL l_{jj} (\cha_j^2-1)} \in \ZZordMih,
\label{1.31}
\qqq
then relations\rx{1.29} are also true
for framed links $\hcL\subset M$.
\end{remark}
\begin{remark}\rm
Similarly to the
TCC invariant
$\ZtrM$ which we mentioned in Introduction,
the definition of $\ZtruaMhL$ was motivated by an asymptotic
approach to the study of the WRT invariant.
The calculation of Witten's path integral for $\ZuaMhL$
in the stationary phase
approximation at large values of $\chK$ suggests that $\ZuaMhL$ splits
into a sum of contributions of connected components $\cM_c$ of the
moduli space $\cM$ of flat $SU(2)$ connections on $M$
\qq
\ZuaMhL = \sum_c \ZcuaMhL,\qquad
\ZcuaMhL = e^{(i\chK/2\pi)\Scsc} \XcuaMhL.
\label{1.32}
\qqq
Here $\Scsc$ is the Chern-Simons invariant of connections of $\cM_c$
and $\XcuaMhL$ are asymptotic power series in $\chKi$ (possibly with
a fractional power of $K$ as a prefactor) whose coefficients depend
on $M,\hcL,c$ and $\uchal$. If $M$ is a \RHS, then the trivial
connection is a separate point in $\cM$ and we conjectured in\cx{Ro1}
that $\ZtruchaMhL$ represents its contribution to the whole WRT
invariant $\ZuaMhL$. \end{remark}

We proved in\cx{Ro1} and\cx{Ro4} that for an empty link $\hcL$,
\qq
\Pebas{0} = 1/\ordhM, \qquad\Pebas{1} = 3\lcwM/\ordhM.
\label{1.33}
\qqq
Then we showed in\cx{Ro5} that an apparent similarity between
\eex{1.14} and\rx{1.33} can be extended.
\begin{theorem}
\label{t1.8}
For an empty link $\hcL$, the coefficients $\Pebas{n}$ of \ex{1.27}
coincide with Ohtsuki's invariants $\l_n(M)$, or in other words,
\qq
\snzi \l_n(M)\,h^n = \ordhMh \,\ZtrM.
\label{1.34}
\qqq
\end{theorem}


\vspace{15pt}

In this paper we will prove the
Lawrence Conjecture\rw{t1.5} extended to the $SO(3)$ WRT
invariant $Z\p$
of links in \RHS\ The definition of this invariant is given by
the following
\begin{theorem}[\cf\cx{KiMe}]
\label{t1.9}
Let $\hcL$ be a framed $L$-component link in a 3-manifold $M$.
Suppose that $M$ is constructed by Dehn's surgery on a
framed $L\p$-component link $\hcLp\subset S^3$. Then for
an odd integer $\chK$ and for a set of odd colors $\uchal$, the
following expression
%
%
\qq
\ZpuchaMhL & = & \imLp\, \chKmLph\,
\emfsLp\,\qtfssLp
\label{1.36}\\
&&\qquad\times
\szubKo \prb{\qpmh{\uuchbe}}\,\JchuaubhLhLp
\nonumber
\qqq
does not depend on the choice of the surgery link $\hcLp$ which is
used to construct $M$.
We call $\ZpuchaMhL$ the $SO(3)$ invariant of $\hcL$ and $M$. It
satisfies the analog of \ex{1.7}
\qq
\begin{array}{cl}
\ZuaMhL = \ZMt\,\ZpuchaMhL & \mbox{if $\eqmod{\chK}{-1}{4}$} \\
\ZuaMhL = \ZMtc\,\ZpuchaMhL & \mbox{if $\eqmod{\chK}{1}{4}$}
\end{array}
\label{1.35}
\qqq
\end{theorem}
\proof
The independence of $\ZpuaMhL$ on the choice of a surgery link
$\hcLp$ can be proved in
exactly the same way as the invariance of $\ZpM$ in\cx{KiMe}. We will
prove \ex{1.35} in Section\rw{s3}. \qed

\def\themt{}
\begin{mt}
\nonumber
For a \RHS\ $M$ let $\chK$ be an odd prime which does not divide
$\ordhM$. If $\hcL$ is a framed link in $M$, then
\qq
\ZpuchaMhL \in \ZZcq
\label{1.40}
\qqq
and
\qq
\ZpuchaMhL \Keq \LegschM
\wdg{\lrbc{ \ordhMh\, \ZtruchaMhLm }},
\label{1.41}
\qqq
or equivalently,
\qq
\ZpuchaMhL \Keq \LegschM \snzi
\PwuchanhL\,\chh^n,
\label{1.42}
\qqq
where $\PuanhL$ are the coefficients in
the expansion\rx{1.27} of
$\ZtruaMhL$.
\end{mt}
\begin{remark}
\label{t1.11}
\rm
Since \eex{1.20} and\rx{1.42} imply that $\l_n(M)=\Pebas{n}$, then
relation\rx{1.19} follows from\rx{1.29}. In fact,\rx{1.29} is a
generalization of\rx{1.19} for the case of a link $\hcL\subset M$
(see also eqs.~(1.87),~(1.90) and~(1.93) in
Theorem~1.8 of\cx{Ro11} for a slightly stronger statement
about the denominators of $\PuanhL$).
\end{remark}

\noindent
{\em Sketch of the proof of Main Theorem.}
Suppose that a \RHS\ $M$ is constructed by a surgery on a framed
algebraically split link $\hcLp\subset S^3$ (\ie all linking numbers
between the components of $\cL\p$ are zero). We prove the theorem
by induction in the number of components of $\hcLp$. An
easy Proposition\rw{t3.3} establishes Main Theorem for $M=S^3$. It
remains to prove that if Main Theorem is true for a \RHS\ $M$, then
it is also true for a \RHS\ $M\p$ constructed by a surgery on a
homologically trivial framed knot $\hcK\subset M$. The key to the
proof is Corollary\rw{t2.2} which describes the structure of
$\ZtruabMhLhK$ and Corollary\rw{2.6} which establishes the similarity
between the gaussian sum of the surgery formula\rx{1.36} and the
gaussian integral of the surgery formula\rx{1.25} as they appear
after the substitution\rx{2.6}.

Our proof of Main Theorem relies exclusively on the properties of the
TCC invariant. We use neither the Murakami-Ohtsuki Theorem\rw{tMO},
nor Theorem\rw{t1.8}. In fact, Main Theorem is stronger than
these theorems (except \ex{1.14} which we will prove in Appendix).

\nsection{Preliminary results}
\label{s2}
All our proofs are based on a particular case of Theorem~1.8
of\cx{Ro11}.
\begin{theorem}
\label{t2.1}
Let $\hcL$ be a framed link and let $\hcK$ be a framed knot with
self-linking number $p$ in a \RHS\ $M$. Let $o$ be
the order of $\cK$ as an element of $\homz$.
For odd colors $\uchal$ of
$\hcL$, there exist the polynomials
\qq
\PuchanMhLKt \in \ZZtohord
\label{2.1}
\qqq
such that
\qq
\PuchanMhLKv{t^{-1}} = \PuchanMhLKt
\label{2.2}
\qqq
and
\qq
\ZtruchabMhLhK =
\ordhMmh\,q^{(1/4)p(\b^2-1)}\,
{q^{\b/2o} - q^{-\b/2o}\over \qpmh{1} }
\snzi {\PuchanMhLKqb\over \DAtnMKqb}\,h^n.
\label{2.3}
\qqq
Here
$\DAMKt$
is the Alexander polynomial of $\cK$ normalized in such a way that
\qq
&\DAMKv{t^{-1}} = \DAMKt,
\label{2.4}\\
&\left. \DAMKt \right|_{t=1} = \ordhM /o.
\label{2.5}
\qqq
\end{theorem}
\proof
For the case of 0-framed link $\cL$ and knot $\cK$, this theorem is
a special case of Theorem~1.8 and Remark~1.10 of\cx{Ro11}.
Remark\rw{t1.7} allows us to extend it to framed links and knots.

\begin{corollary}
\label{t2.2}
If $\hcK$ is a framed homologically trivial knot with self-linking
number $p$ in a \RHS\ $M$, then for odd colors $\uchal$
\qq
\ZtruchabMhLhK & = &
\ordhMmh\,q^{1/2}\,h^{-1}\,q^{(1/4)p(\b^2-1)}
\label{2.6}\\
&&\qquad\times
\smnz \duchamnMhLK\,(\qpmh{\b})^{2m+1}\,h^n
\nonumber\\
\lefteqn{
\duchamnMhLK \in \ZZhordMi
}
\label{2.7}
\qqq

\end{corollary}

\proof
Since $\cK$ is homologically trivial, then $o=1$ and also
\qq
\DAMKt \in \ZZtpm.
\label{2.8}
\qqq
Therefore relations\rx{2.2},\rx{2.4} and\rx{2.5} imply that $\DAMKt$
and $\PuchanMhLKt$ are polynomials of $(\thpm)^2$
\qq
&\DAMKt = \ordhM + \smo a_m\,(\thpm)^{2m},\qquad a_m\in \ZZ,
\label{2.9}\\
& \PuchanMhLKt = \smz b_m\,(\thpm)^{2m},\qquad
b_m\in\ZZhordMi.
\label{2.10}
\qqq
The expansion of denominators in the sum of \ex{2.3} in powers of
$(\qpmh{\b})^2$ leads to \ex{2.6}.\qed

To prove Main Theorem, we will also need some simple facts.
\begin{lemma}[H.~Murakami\cx{Mu1}]
\label{t2.3}
If $f(q)\in \ZZqpm$ and its expansion in powers of $h=q-1$ is
\qq
f(q) = \snzi a_n\,h^n,
\label{2.11}
\qqq
then for any $N>0$
\qq
f(q) = \snzN a_n\,h^n + h^{N+1}x, \qquad
\mbox{where $x\in \ZZqpm.$}
\label{2.12}
\qqq
\end{lemma}
\begin{lemma}
\label{t2.4}
If $\chK$ is prime and $m,n\in\ZZ$, $n$ not divisible by $\chK$, then
there is an equlity in $\ZZKcq$
%
\qq
\chq^{mn\s} \Keq (q^{m/n})\wve,
\label{2.13}
\qqq
where $q^{m/n}$ is understood as a power series in $h$
\qq
& q^{m/n} = (1+h)^{m/n} = \skzi {m/n \choose k} h^k
\label{v1.1}
\\
&
{m/n \choose k} ={ {\fmn} \lrbc{\fmn - 1}\cdots \lrbc{\fmn - k +
1} \over k! } =
{ m(m-n)\cdots(m-n(k-1)) \over n^k\, k!}.
\label{v1.2}
\qqq
\end{lemma}
\proof
It is easy to see that if $n$ is not divisible by $\chK$, then
\qq
{m/n \choose k} \in \QchK
\label{v1.3}
\qqq
and as a result
\qq
{m/n \choose k}\wve = {m\,(1/n)\wve \choose k}\in\ZZK.
\label{v1.4}
\qqq
Therefore according to \ex{v1.1},
the \rhs of \ex{2.13} can be presented in $\ZZKcq$ as a convergent
power series in $\chh$
\qq
(q^{m/n})\wve = \skzi {m/n \choose k}\wve \chh^k =
\skzi {m\,(1/n)\wve \choose k} \chh^k.
\label{v1.4*}
\qqq

Since $\chq^{m n\s}$ does not depend on the choice of $n\s$
which satisfies \ex{1.37}, then according to definition\rx{1.12}, for
any positive integer $N$ we can set
$n\s = \nssN$
and present the \lhs of \ex{2.13} in $\ZZKcq$ also as a
convergent power series in $h$
%
%
\qq
\chq^{m n\s} = \chq^{m\, \nssN} = (1+\chh)^{m\,\nssN}
=\skzi {m\,\nssN \choose k} \chh^k.
\label{v1.5}
\qqq
Now \ex{2.13} follows from the fact that each term of the
series\rx{v1.5} converges $\chK$-adicly to the corresponding term of
the series\rx{v1.4*} as $N\longrightarrow \infty$. \qed

\def\epra#1{ #1_{\rm asympt.} }
\def\eprk#1{ #1_{\rm cycl.} }

\def\pmchK{ (p,m,\chK) }
\def\pmK{ (p,m,K) }

\def\pmchq{ (p,m,\chq) }
\def\pmq{ (p,m,q) }

\def\Yc{ \eprk{Y} }
\def\Ycv{ \Yc\pmchq }
\def\Ya{ \epra{Y} }
\def\Yav{ \Ya\pmq }

Our proof of the Lawrence conjecture is based on a simple relation
between gaussian sums and gaussian integrals.
We introduce the following notation: for a function
$f(\chb)$,
%
\qq
\spzbKo f(\chb) = \szbKo f(\chb) + {1\over 2}\,f(\chK).
\qqq
%

\begin{lemma}
\label{t2.5}
For $\chK$ being an odd prime and for $p,m\in \ZZ$, $p$ not
divisible by $\chK$, consider the expression
\qq
&\eprk{X}  =
\chKmh\,\emfsp
\spzbKo \qfspb \chqmbp\in\IC
\label{vr1}
\qqq
and a function of $q$
\qq
&\epra{X}  =
\eKmoph\,\pah\,\empfsp \intmzz \qfpb\qmbp d\b.
\label{vr2}
\qqq
which is well-defined for $0<|q|<1$.
%
%
Then
\qq
&
\eprk{X} =  \Legschp \, \qmmps \in \ZZcq\subset \IC,
\label{2.14}
\qqq
while\rx{vr2} can be extended analyticly to the vicinity
of $q=1$
and expanded in power series of $h=q-1$
\qq
&
\epra{X} = q^{-m^2/p} \in \ZZhpi.
\label{2.14*}
\qqq
%
Finally,
\qq
\eprk{X} = \Legschp \lrbc{\epra{X}}\wve.
\label{2.15}
\qqq
\end{lemma}

\proof
We calculate the sum
\qq
\spzbKo \qfspb \chqmbp
& = &
\szbmK \chq^{4\s p\chb^2 + m\chb}
= \szbmK \chq^{4\s p\,(\chb+2mp\s)^2 - m^2 p\s}
\nonumber\\
& = & \qmmps \szbmK \qfspb
= \qmmps \sgmKh \chq^{4\s p\,(2\chgamma+1)^2}
\nonumber\\
& = &
\qmmps \sgmKh \chq^{p\,(\chgamma+2\s)^2}
= \qmmps \sgoK \chq^{p\chgamma^2}
\nonumber\\
& = &
\Legschp\,e^{(i\pi/4)(1-\kappa)\prosign(p)}\,\chKh\,\qmmps
\label{2.16}
\qqq
In the last line we used the formula for the gaussian sum which can
be found, for example, in\cx{Mu5}, Chapter~6:
\qq
\sgoK \chq^{p\chgamma^2}
=\Legschp\,e^{(i\pi/4)(1-\kappa)\prosign(p)}\,\chKh.
\qqq
Equation\rx{2.14}
follows from \ex{2.16}.

Next we calculate the gaussian integral
\qq
\intmzz\qfpb\qmbp\, d\b
& = & 2\qmmp \intmzz q^{(1/4)\,p\,(\b+2m/p)^2}\,d\b
\nonumber\\
& = & 2\qmmp \intmzz e^{(i\pi/2K)p\b^2}d\b
\nonumber\\
& = & (8K)^{1/2}\,e^{(i\pi/4)\prosign(p)}\,|p|^{-1/2}\,\qmmp.
\label{2.17}
\qqq
Relation\rx{2.14*} follows from this equation, because
$q^{-m^2/p}=(1+h)^{-m^2/p}\in\ZZhpi$.
Equation\rx{2.15} follows from \eex{2.14},\rx{2.14*} and from
Lemma\rw{t2.4}, because the \rhs of \ex{2.15} is equal to
$\Legschp\lrbc{ \qmmp }\wve$. \qed
\begin{corollary}
\label{t2.6}
For $\chK$
being an odd prime, $m,p\in \ZZ$, $p$ not divisible by
$\chK$, consider the expression
\qq
\Ycv & = &
\chhi\,\chKmh\,\emfsp \szbKo \qfspb\chqmbtm\in\IC
\label{vr3}
\qqq
and a function
\qq
\Yav & = &
\hi\,\eKmoph\,\pah\,\empfsp
\nonumber\\
&&\qquad\times\intmzz\qfpb\qmbtm\, d\b.
\label{vr4}
\qqq
which is well-defined for $0<|q|<1$. Then $\Ycv\in\ZZcq\subset\IC$
and it is $\chK$-adicly small
\qq
\Ycv & = &
\chh^m\,x,
\qquad x\in \ZZcq.
\label{2.18}
\qqq
$\Yav$ can be extended analytically to the vicinity of $q=1$
and expanded there in powers of $h=q-1$, so that $\Yav\in\ZZhpi$ and
it is asymptotically small
%
\qq
\Yav & = &
h^m\,x\p,
\qquad
x\p
\in \ZZhpi.
\label{2.18*}
\qqq
Finally,
\qq
\Ycv = \Legschp \lrbc{\Yav}\wve
\label{2.19}
\qqq
\end{corollary}


\proof
Obviously,
\qq
\chq^{\chK/2} - \chq^{-\chK/2} = 0,
\qqq
so we can replace the sum
$\szbKo$  in the expression\rx{vr3} for
$\eprk{Y}$ by the extended sum $\spzbKo$.
Since for $m\geq 0$
\qq
&\qmbtm = \snzmo a_{m,n}\,(q^{n\b} + q^{-n\b})
,\qquad a_{m,n}\in\ZZ,
\label{2.20}
\qqq
then we conclude from \ex{2.14} and\rx{2.15} that
\qq
&
\chh\, \eprk{Y} \in\ZZcq,\qquad
h\, \epra{Y} \in \ZZhpi,\qquad
\chh\, \eprk{Y} = \Legschp \lrbc{ h\,\epra{Y} }\wve
\label{2.21}
\qqq
Since
\qq
\Kmh \intmzz\qfpb\qmbtm\, d\b
= O(K^{-m-1})\quad\mbox{as $K\rightarrow \infty$,}
\label{2.23}
\qqq
then the expansion of $h\,\epra{Y}$ in powers of $h$ starts at
$h^{m+1}$. Therefore $\epra{Y}$ is divisible by $h^m$ in $\ZZhpi$,
$\eprk{Y}$ is divisible by $\chh^m$ in $\ZZcq$ and the
relations\rx{2.18}--(\ref{2.19}) follow from\rx{2.21}.\qed

\nsection{Proof of Main Theorem}
\label{s3}
We start with Symmetry Principle formulated by R.~Kirby and P.~Melvin
in\cx{KiMe}.
\begin{theorem}[R.~Kirby, P.~Melvin\cx{KiMe}]
\label{t3.1}
Let $\hcL$ be a framed $L$-component link in $S^3$. For
a fixed integer $j$, $1\leq j\leq L$,
consider two sets of colors $\uchal$ and $\uchal\p$
such that
\qq
\cha\p_k = \cha_k\;\;\;\mbox{if $k\neq j$}, \qquad
\cha_j\p = \chK-\cha_j.
\label{3.1}
\qqq
Then
\qq
\JuchaphL = i^{\kappa l_{jj}}\,(-1)^{l_{jj}\cha_j +
\sokLnj l_{jk}(\cha_k - 1)}\,\JuchahL,
\label{3.2}
\qqq
here $l_{jk}$ are the linking numbers of $\hcL$ while $\kappa$ is
defined by \ex{1.38}.
\end{theorem}
\begin{remark}
\label{t3.2}
\rm
This theorem follows easily from the Main Theorem of\cx{Ro10} which
is a particular case of Theorem~1.8 of\cx{Ro11} (see
Appendix~C of\cx{Ro10} for the details).
\end{remark}


\pr{\ex{1.35} of Theorem}{t1.9}
Our proof is a slightly shortened version of the proof of \ex{1.7}
in\cx{KiMe}.
The Symmetry Principle applied to the colors $\uuchbe$ in the surgery
formula\rx{1.5} suggests that the sum over $\uuchbe$ can be split into
a sum over the elements of the symmetry group generated by the
transformations\rx{3.1} acting on $\uuchbe$, and over the orbits of
that action. Since according to our assumptions, $\chK$ is odd, then
each orbit contains a set of colors $\uuchbe$ which are all odd.
Therefore for odd colors $\uchal$ we can rewrite \ex{1.5} as
\qq
\ZuaMhL & = & \imLp\,\tchKmLph\,\emtfsLp\,\qchtfsLp
\label{3.3}\\
&&\qquad\times
\XhcLpk \szubKo\prb{\qchpmh{\uuchbe}}\,\JchuaubhLhLp,
\nonumber
\qqq
where
\qq
\XhcLpk & = &
\sum_{J\subset\{1,\ldots,L\p\}}
i^{\kappa \sum_{j\in J} l\p_{jj}}
\,(-1)^{\sum_{j,k\in J,\,j\leq k}l\p_{j,k} }
\label{3.4}
\qqq
and $l\p_{ij}$ are the linking numbers of $\hcLp$.
Note that
\qq
\XhcLp{1} = \overline{\XhcLp{-1}}
\label{3.5}
\qqq
and the sum $\XhcLpk$ depends on $K$ only through $\kappa$.
Consider the invariant $\ZM$ at $K=3$. Since
$\left.\JuchahL\right|_{\uchal=1} = 1$, then according to \ex{3.3},
\qq
\ZMt = 2^{-L\p/2}\,e^{-(i\pi/4)\prosign(\hcLp)}\,\XhcLp{-1}.
\label{3.6}
\qqq
%
A combination of \eex{3.3},\rx{3.5} and\rx{3.6} shows that \ex{1.35}
is satisfied if the invariant $\ZpuchaMhL$ was given by the
expression
\qq
\ZpuchaMhL
& = &
\imLp\,\chKmLph\,e^{-(i\pi/4)(\kappa+3)\prosign(\hcLp)}\,\qchtfsLp
\nonumber
\\
&&\qquad\times
\szubKo \prb{\qchpmh{\uuchbe}} \JchuaubhLhLp.
\label{3.7}
\qqq
In fact, since
\qq
\eqmod{{1-\kappa \chK\over 4}}{4\s}{\chK},
\label{3.8}
\qqq
then it is easy to see that
\qq
e^{-(i\pi/4)(\kappa+3)\prosign(\hcLp)}\,\qchtfsLp
=
\emfsLp\,\qtfssLp
\label{3.9}
\qqq
and \ex{3.7} is equivalent to \ex{1.36}.\qed

\begin{proposition}
\label{t1.10}
Suppose that $\chK$ is an odd
integer. Let $\hcK\subset M$ be a framed
knot with self-linking number $p\neq 0$. If $M\p$ is constructed by
Dehn's surgery on $\hcK$, then for odd colors $\uchal$
\qq
\ZpuchaMphL & = & \chKmh\,\sip\,\emfsp\,\qtfssp
\label{1.39}\\
&&\qquad\times
\szbKo (\qpmh{\chb})\,\ZpuchabMhLhK.
\nonumber
\qqq
\end{proposition}

\proof
Suppose that $M$ is constructed by a surgery on a framed link
$\hcLp\subset S^3$. Then $M\p$ can be constructed by a surgery on
$\hLphK$ with $\hcK\subset S^3$ having a self-linking number
\qq
p_0 = p + \sjkoLp l\p_{0j}\,l\p_{0k}\,(l\p)^{-1}_{jk},
\label{3.10}
\qqq
here $l\p_{0j}$ are the linking numbers between $\cK$ and components
of $\hcLp$ while $(l\p)^{-1}_{jk}$ is the inverse linking matrix of
$\hcLp$. Since
\qq
\prosign(\hLphK) = \prosign(p) + \prosign(\hcLp)
\label{3.11}
\qqq
and for $p\neq 0$
\qq
i^{-1}\,
\emfs{p}{+} = \prosign(p)\,\emfsp,
\label{3.12*}
\qqq
then comparing \ex{1.36} for the surgeries on $\hcLp$ and on $\hLphK$
we arrive at \ex{1.39}.\qed

The following two propositions will allow us to prove Main Theorem by
induction.
\begin{proposition}
\label{t3.3}
The Main Theorem is true for $M=S^3$.
\end{proposition}

\proof
According to \ex{1.36} applied in the trivial case of an empty surgery
link,
\qq
\Zpbasch{\uchal}{S^3,\hcL} = \JuchahL.
\label{3.12}
\qqq
Therefore relation\rx{1.40} follows from\rx{1.2}. Equation\rx{1.22}
indicates that
$\Ztrbas{\uchal}{S^3,\hcL}$ is the expansion of $\JuahL$
in powers of $h$, so \ex{1.41} follows from Lemma\rw{t2.3} together
with \ex{v2.4}
and the definition of $\chK$-adic convergence\rx{1.18}.\qed
\begin{proposition}
\label{t3.4}
Suppose that a \RHS\ $M\p$ is constructed by a surgery on a
homologically trivial framed knot $\hcK$
 in another \RHS\ $M$. If the Main
Theorem is true for $M$, then it is true for $M\p$.
\end{proposition}
\proof
Suppose that $\ordhMp$ is not divisible by an odd prime $\chK$.
Since
\qq
\ordhMp = |p|\,\ordhM,
\label{sp.3}
\qqq
where $p$ is the self-linking number of $\hcK$,
then neither $p$ nor $\ordhM$ is divisible by $\chK$.
We assume that \ex{1.42} is true for the link $\hLhK\subset M$
\qq
\ZpuchabMhLhK \Keq \LegschM\snzi \Pwbas{\uchal,\chb\sepx n}{\hLhK}\,
\chh^n.
\label{3.13}
\qqq
On the other hand, if $\uchal$ and $\chb$ are odd, then in view of
\eex{1.27},\rx{2.6} and
Lemma\rw{2.3}, for any $N>0$
\qq
\lefteqn{
\snzN \Pbas{\uchal,\chb\sepx n}{\hLhK}\,h^n}
\label{3.14}\\
& = &
\chq^{1/2}\,
\chhi\,\chq^{(1/4)p\,(\chb^2-1)}\duchamnMhLK\,\smnzN\chqmbto\,\chh^n
+ \chh^{N+1}\xon,
\nonumber
\qqq
where
\qq
&\duchamnMhLK\in\ZZhordMi,\qquad
\xon\in \ZZqhord.
\label{3.15}
\qqq
Since we assumed that $\chK$ does not divide $\ordhM$, then the first
relation of\rx{3.15} implies that $\duchamnMhLK\in\QchK$. Therefore
$\dchv$ is well-defined and we can also use $\remkN{\dchv}$
as its $\chK$-adic approximation. Thus
equations\rx{3.13} and\rx{3.14} together with the relations\rx{v2.4}
and
\qq
\chq^{(1/4)p\,(\chb^2-1)}
= \chq^{4\s p\,(\chb^2-1)}\;\;\mbox{if $\chb$ is
odd},\qquad \chq^{1/2} = -\chq^{2\s}
\label{3.16}
\qqq
imply that for any $N_0$ there exists $N_0\p$ such that for any
$N>N_0\p$
\qq
\lefteqn{
\ZpuchabMhLhK
=
-\LegschM\,\chq^{2\s}\,\chhi\,\chq^{4\s p\,(\chb^2-1)}
}
\label{3.17}\\
&&\times
\smnzN
\remkN{\dchv}\,
\chqmbto\,\chh^n +
\chK^{N_0+2}\xtw,\qquad \xtw\in\ZZcq.
\nonumber
\qqq
%
We substitute
\ex{3.17}
into \ex{1.39}. Since
\qq
\chKh \,e^{(i\pi/4)(\kappa-1)} =
\sgoK \chq^{\chgamma^2}\in \ZZcq
\label{3.18}
\qqq
and for odd $\chb$
\qq
\qchpmh{\chb} = \chq^{1/2}(\chq^{(\chb-1)/2} - \chq^{-(\chb+1)/2}) =
- \chq^{2\s} (\chq^{(\chb-1)/2} - \chq^{-(\chb+1)/2})\in \ZZcq,
\label{3.19}
\qqq
then it is easy to see that the contribution of the term
$\chK^{N_0+2}\,\xtw$ to $\ZpuchaMphL$ is equal to
\qq
\lefteqn{
\chKmh\,\sip\,\emfsp\,\qtfssp
\szbKo (\qchpmh{\chb})\,\chK^{N_0+2}\,\xtw}\hspace{4in}
\label{3.20}\\
&=&
\chK^{N_0+1}\,\xth,
\qquad
\xth\in \ZZcq.
\nonumber
\qqq
Thus we come to the following statement: for any $N_0>0$ there exists
$N$ such that
\def\qassts{ \chq^{2\s-4\s p + 3\cdot 4\s\prosign(p)} }
\def\qassts{ \chq^{4\s(2- p + 3\prosign(p))} }
\qq
\lefteqn{\ZpuchaMphL =}&&
\nonumber\\
& = & - \LegschMp \prosign(p)\,
\qassts
\smnzN
\remkN{\dchv} \Ycv\,\chh^n
\nonumber\\
&&\qquad+ \chK^{N_0+1}\,\xth,\qquad \xth\in\ZZcq,
\label{3.23}
\qqq
where $\Ycv$ is defined by \ex{vr3}. Since \ex{2.18} implies that
$\Ycv\in\ZZcq$, then the claim\rx{1.40} of Main Theorem for the
manifold $M\p$ follows easily from \ex{3.23}.

It remains to prove \ex{1.41} for $M\p$. If we substitute the
expression for $\ZtruabMhLhK$ from \ex{2.6}
into the surgery
formula\rx{1.25} and use the relations\rx{sp.3} and
\qq
i\,e^{-(i\pi/2)\sip} = \sip\;\;\mbox{if $p\neq 0$},
\label{3.21}
\qqq
then we find that
\def\qasst{ q^{(2 - p + 3\prosign(p))/4} }
\qq
\ordhMph\,\ZtruchaMphL & = & -
\prosign(p)\, \qasst
\hspace{-7pt}\smnzN \hspace{-8pt}\duchamnMhLK\,\Yav\,h^n
\nonumber\\
&&\qquad\qquad+
h^{N+1}\,\xfo,\qquad
\xfo\in\ZZhordMpi[[h]],
\label{3.23*}
\qqq
where $\Yav$ is defined by \ex{vr4}.
The remainder $h^{N+1}\xfo$ in this formula comes from the
terms $\duchamnMhLK\,\Yav$ with $m+n>N$
which were not included in the
sum of \ex{3.23*}. Their contribution is estimated with the help of
relation\rx{2.18*}.

Now if we compare \eex{3.23}
and\rx{3.23*}, then \ex{1.41} for $M\p$ follows easily from \ex{2.19}
if we recall
\ex{1.12} which relates $\dchv$ and $\remkN{\dchv}$ and if we use the
relation
\qq
\lrbc{\qasst}\wve = \qassts
\qqq
which follows from Lemma\rw{t2.4}. \qed

\begin{corollary}
\label{t3.5}
Main Theorem is true for a \RHS\ $M$ which is constructed by Dehn's
surgery on an algebraically split framed link $\hcLp\subset S^3$.
\end{corollary}
\proof
The proof is by induction in the number $L\p$ of components of
$\hcLp$. Proposition\rw{t3.3} demonstrates that the claim is true for
$L\p=0$, that is, when $\hcLp$ is an empty link and as a result
$M=S^3$. Suppose that the claim is true for
$L\p$-component links. Consider an $(L\p+1)$-component algebraically
split link $\hcLp$. The surgery on the first $L\p$ components of
$\hcLp$ produces a \RHS\ $M\p$,
and the $(L\p+1)$-st component of $\hcLp$ is
homologically trivial inside $M\p$. Therefore Main Theorem is
true
for the surgery on the whole link $\hcLp$ because of
Proposition\rw{t3.4}.\qed
\begin{lemma}
\label{t3.6}
Suppose that a \RHS\ $M\p=M\#\Lpo$ is a connected sum of a lens space
$\Lpo$ and of a \RHS\ $M$ which contains a framed link $\hcL$.
%
%
Suppose that Main Theorem holds for $M\p$. If an odd prime $K$
divides neither $p$, nor $\ordhM$,
then the claim of Main Theorem is true for $M$.
\end{lemma}
\proof
The lens space $\Lpo$ can be constructed by Dehn's surgery on a
framed unknot in $S^3$ with self-linking number $-p$. Since
\qq
J_{\chb}({\rm unknot};q) = \chq^{-p\,(\b^2-1)/4}\,
{\qpmh{\chb}\over \qpmh{1}},
\label{3.25}
\qqq
then according to \ex{1.25}
\qq
\ZtrLp
& = &
-i\,\eKmoph\,e^{(3\pi i/4)\sip}\,
q^{-(3/4)\sip}\,q^{p/4}\,
(\qpmh{1})^{-1}
\nonumber\\
&&\qquad\times\intmzz e^{-(i\pi/2K)p\b^2}\,
(\epmh{\b})^2\,d\b
\nonumber\\
& = & \sip\, |p|^{-1/2} \,q^{(p+2/p-3\sip)/4}\,
{q^{1/2p} - q^{-1/2p} \over \qpmh{1}}
\label{3.26}
\qqq
(this formula was first derived by L.~Jeffrey in\cx{Je}).
Corollary\rw{t3.5} indicates that Main Theorem is true for $\Lpo$, so
since
$
\ordhozbas{\Lpo} = |p|
$
then
\qq
\ZpLp & = & \Legschp \lrbc{
\pah\,\ZtrLp}\wve
\nonumber\\
& = &
\Legschp\,\sip\,\chq^{4\s(p+2p\s-3\sip)}\,
{\chq^{2\s p\s} - \chq^{-2\s p\s} \over \chq^{2\s} - \chq^{-2\s} }.
\label{3.28}
\qqq
Note that since
$
\chq^{2\s} = (\chq^{2\s p\s})^p,
$
then
\qq
{\chq^{2\s} - \chq^{-2\s} \over
\chq^{2\s p\s} - \chq^{-2\s p\s}} \in \ZZcq
\label{3.30}
\qqq
and as a result
\qq
{1\over\ZpLp} \in \ZZcq.
\label{3.31}
\qqq

Suppose that an odd prime $K$ divides neither $p$, nor $\ordhM$.
Since
$\ordhozbas{\Lpo} = |p|$, then in view of multiplicativity
of $\ordhM$ under the operation of connected sum,
\qq
\ordhMp = \ordhozbas{\Lpo}\,\ordhM = |p|\,\ordhM
\label{3.36}
\qqq
and, as a result, $K$ does not divide $\ordhMp$. Therefore, according
to our assumptions, the claims of Main Theorem are true for $M\p$
%
\qq
&\ZpuchaMphL \in \ZZcq,
\label{3.32}\\
&\ZpuchaMphL = \LegschMp\lrbc{ \ordhMp^{1/2}\,\ZtruchaMphL}\wve.
\label{3.33}
\qqq
Both quantum invariants
$\ZpuaMhL$ and $\ZtruaMhL$ are also multiplicative under the operation
of connected sum, so
\qq
&\ZpuchaMhL = \ZpuchaMphL/\ZpLp,
\label{3.34}\\
&\ZtruaMhL = \ZtruaMphL/\ZtrLp,
\label{3.35}
\qqq
Therefore relation\rx{1.40} follows from \ex{3.34} and from
relations\rx{3.32} and\rx{3.31}, while \ex{1.41} follows from
\eex{3.35},\rx{3.33} and\rx{3.28}.\qed

The last statement that we need in order to prove Main Theorem is the
lemma that H.~Murakami introduced in\cx{Mu2} at the suggestion of
T.~Ohtsuki.
\begin{lemma}[H.~Murakami\cx{Mu2}]
\label{tMu}
For a \RHS\ $M$ and an odd prime $\chK$ there
exists an algebraically
split framed link $\hcLp\subset S^3$ such that Dehn's surgery on
$\hcLp$ produces a \RHS\ $M\p$ which is a connected sum of $M$ and of
a finite number of lens spaces $L_{p_j,1}$, $1\leq j\leq N$,
\qq
M\p = M\#L_{p_1,1}\#\cdots\#L_{p_N,1},
\label{3.37}
\qqq
such that neither of $p_j$ is divisible by $\chK$.
\end{lemma}
\noindent{\em Proof of Main Theorem.\/}
For a \RHS\ $M$ let
$\chK$ be an odd prime number which does not divide
$\ordhM$. Consider the \RHS\ $M\p$ of Lemma\rw{tMu} which can be
constructed by Dehn's surgery on a framed algebraically split link
$\hcLp\subset S^3$. Corollary\rw{t3.5} says
that Main Theorem is true for $M\p$ and then Lemma\rw{t3.6} implies
that it is also true for $M$.\qed

\noindent
{\bf Acknowledgements.}
I am thankful to R.~Lawrence and A.~Vaintrob
for their comments and suggestions. I am greatly indebted to
A.~Vaintrob for his persistence and patience in helping me to improve
the exposition of this paper. This work was supported by NSF Grant
DMS-9704893.


\section*{Appendix}
\appendix
\setcounter{equation}{0}
\section{Casson's invariant as the second coefficient of the TCC
invariant}

In\cx{Ro1} we used path integral arguments in the framework of the
Chern-Simons quantum field theory of\cx{Wi} in order to derive a much
weaker version of Theorem\rw{t2.1}. Considering that result as a
conjecture, we used it to derive relations\rx{1.33}. Now that we have
Theorem\rw{t2.1} at our disposal, we repeat the arguments of\cx{Ro1}
in order to present a rigorous proof of\rx{1.33} which in conjunction
with Theorem\rw{t1.8} constitutes an alternative proof of Murakami's
result\rx{1.14}.


The first of \eex{1.33} was already established in\cx{Ro11}, so it
remains to show that for a \RHS\ $M$
\qq
{1\over 3}\,\ordhM\,\Pebas{1} = \lcwM.
\label{a2.1}
\qqq

We have to establish some preliminary facts. We proved in\cx{Ro11}
that for a knot $\cK\subset M$
\qq
\Pbas{0}{\cK;t} = 1,
\label{A2.2}
\qqq
here $\Pbas{n}{\cK;t} = 1$ are the polynomials of \ex{2.3} in the
case of an empty link $\hcL$. A simple corollary of \eex{A2.2}
and\rx{2.5} is that if $\cK$ is homologically trivial in $M$, then
\qq
&
\dzz = 1/\ordhM,
\label{a2.3}\\
&
\doz = {1\over 2}\,\ordhMmt\,\DAddMK.
\label{a2.4}
\qqq
Here
\qq
\DAddMK = {d^2\over d t^2} \left. \DAMKt\right|_{t=1}
\label{a2.5}
\qqq
and the coefficients $\dMKbas{m,n}$ come from \ex{2.6} written for an
empty link $\hcL$.

Suppose that a \RHS\ $M\p$ is constructed by a surgery on a framed
homologically trivial knot $\hcK\subset M$ with self-linking number
$p$. Then according to \ex{1.25},\rx{1.27} and\rx{2.6}
\qq
\lefteqn{
\PeMpbas{0} + h\,\PeMpbas{1} =
-i\, \eKmoph\,\emtfsp\,\qtfsp\,\pah\,(\qpmh{1})^{-1}
}
\nonumber\\
&&\times\intmzz q^{(1/4)p\,(\b^2-1)}\lrbc{
\dzz\,(\qpmh{\b})^2 + \doz\,(\qpmh{\b})^4
\right.
\nonumber\\
&&\left.
\hspace{2in}
+
h\,\dzo\,(\qpmh{\b})^2 }d\b + O(h^2)
\label{a2.6}
\qqq
and as a result,
\qq
\PeMpbas{1} = {1\over |p|}
\lrbs{
\lrbc{ {3\over 4}\,\sip - {p\over 4} - {1\over 2p} }\dzz
+
\dzo - {6\over p}\,\doz },
\label{a2.7}
\qqq
or in view of \eex{a2.3},\rx{a2.4}
\qq
\hspace{-20pt}
\PeMpbas{1} = {1\over\ordhMp}
\lrbc{ {3\over p}\,\sip - {p\over 4} - {1\over 2p} -
{3\over 4}\,{\DAddMK\over \ordhM} + \ordhM\,\dzo}.
\label{a2.8}
\qqq
%
\begin{proposition}
\label{at1}
Equation\rx{a2.1} is true for a \RHS\ $M$
constructed by a surgery on a
framed algebraically split link $\hcLp\subset S^3$.
\end{proposition}
\proof(\cf the proof of Corollary\rw{t3.5})
We prove the proposition
by induction in the number of components $L\p$ of $\hcLp$.

If $L\p=0$, then $M=S^3$. Since
\qq
\lcw(S^3) =0,
\label{a2.9}
\qqq
then \ex{a2.1} is true in view of \eex{sp.1} and\rx{1.22}.

Suppose that \ex{a2.1} is true for all \RHS\ constructed by surgeries
on $L\p$-component links in $S^3$. Let us prove that it is also true
for manifolds constructed by surgeries on $(L\p+1)$ component links.
It is enough to prove the following: if \ex{a2.1} is true for a
\RHS\ $M$ then it is also true for a \RHS\ $M\p$ constructed by
Dehn's surgery on a framed homologically trivial knot
$\hcK\subset M$.

Equations\rx{1.27},\rx{1.28} imply that
\qq
\Ztrbas{\b}{M,\hcK} = \ordhMmh
\sum_{m,n\geq 0\atop |\um|\leq n/2} \dmnK\,\b^{2m+1}\,h^{n}.
\qqq
If we substitute $\b=1$ into this equation and take into account
the relation
%
\qq
&\left.\Ztrbas{\b}{M,\hcK} \right|_{\b=1} = \ZtrM,
\label{a2.10}
\qqq
%
then we conclude that
%
\qq
\dzo = \Pebas{1}. 
\label{a2.12}
\qqq

Assume that \ex{a2.1} is true for $M$. Then in view of \ex{a2.12}
\qq
\dzo = 3\lcwM/\ordhM,
\qqq
and therefore
according to \ex{a2.8},
\qq
(1/3)\,\ordhMp\,\PeMpbas{1} =
{1\over 4}\,\sip  - {p\over 12} - {1\over 6p} -
{1\over p}\,{\DAddMK\over\ordhM} + \lcwM.
\label{a2.13}
\qqq
According to\cx{Wa}, the \rhs of this equation is equal to
$\lcw(M\p)$, hence we proved \ex{a2.1} for $M\p$.\qed

Since a lens space $\Lpo$ is constructed by Dehn's surgery on an
unknot in $S^3$, then Proposition\rw{at1} implies that \ex{a2.1} is
true for $\Lpo$. Suppose that \ex{a2.1} holds for a manifold
%
$M\p = M\#\Lpo$,
%
with $M$ being a \RHS\ Since both $\lcw(M)$ and
$(1/3)\,\ordhM\,\PeMpbas{1}$ are additive under the operation of
connected sum (the latter additivity follows from the
multiplicativity of $\ZtrM$), then \ex{a2.1} should also be true for
$M$. In view of Lemma\rw{tMu} and Proposition\rw{at1}, this proves
\ex{a2.1} for any \RHS\ (\cf Proof of Main Theorem).

\end{document}